\documentclass{aa}
\usepackage[varg]{txfonts}

\usepackage{natbib}
\bibpunct{(}{)}{;}{a}{}{,}

\usepackage[utf8]{inputenc}
\usepackage{amsmath,amsfonts}
\usepackage{algorithm}
\usepackage{algpseudocode}

\usepackage{graphicx}
\graphicspath{ {figures/} }

\usepackage{setspace}
\numberwithin{equation}{section}

\newcommand{\signal}[1]{\left\lceil #1 \right\rfloor}
\newcommand{\npix}{{\ensuremath{n_{\mathrm{pix}}}}}
\newcommand{\nfreq}{{\ensuremath{n_{\mathrm{freq}}}}}
\newcommand{\nt}{{\ensuremath{n_t}}}

\usepackage{bbm}
\newcommand{\sig}[1]{\mathsf{#1}}
\newcommand{\sigg}[1]{\mathbbmss{#1}}
\newcommand{\A}{\sigg{A}}
\newcommand{\M}{\sigg{B}}

\usepackage{hyperref}
\usepackage{color}

\definecolor{linkcolour}{rgb}{0.64,0,0}
\hypersetup{colorlinks,breaklinks,urlcolor=linkcolour, linkcolor=linkcolour}
\hypersetup{draft}

\usepackage{cleveref}

\date{April 29, 2020} 

\begin{document}

\titlerunning{Accelerating linear system solvers for time-domain component separation of CMB data}
\authorrunning{J. Pape\v{z}, L. Grigori, R. Stompor}

\title{
Accelerating linear system solvers for time-domain component separation of cosmic microwave background data
}

\author{J. Pape\v{z}\inst{1}\thanks{\email{jan@papez.org}}, L. Grigori\inst{2}, R. Stompor\inst{3,4}}
\institute{
INRIA Paris, Sorbonne Universit\'e, Universit\'e Paris-Diderot SPC, CNRS, Laboratoire Jacques-Louis Lions, ALPINES team, France\\ Currently at
Institute of Mathematics, Czech Academy of Sciences, Prague, Czech Republic \and
INRIA Paris, Sorbonne Universit\'e, Universit\'e Paris-Diderot SPC, CNRS, Laboratoire Jacques-Louis Lions, ALPINES team, France
\and Universit\'e de Paris, CNRS, AstroParticule et Cosmologie, F-75013 Paris, France
\and
CNRS-UCB International Research Laboratory, "Centre Pierre Bin\'etruy", UMI2007, CPB-IN2P3
}

\abstract{
Component separation is one of the key stages of any modern cosmic microwave background (CMB) data analysis pipeline. It is an inherently nonlinear procedure and typically involves a series of sequential solutions of linear systems with similar but not identical system matrices, derived for different data models of the same data set. Sequences of this type arise, for instance, in the maximization of the data likelihood with respect to foreground parameters or sampling of their posterior distribution. However, they are also common in many other contexts. 
In this work we consider solving the component separation problem directly in the measurement (time-) domain. This can have a number of important benefits over the more standard pixel-based methods, in particular if non-negligible time-domain noise correlations are present, as is commonly the case.
The approach based on the time-domain, however, implies significant computational effort because the full volume of the time-domain data set needs to be manipulated. To address this challenge, we propose and study efficient solvers adapted to solving time-domain-based component separation systems and their sequences, and which are capable of capitalizing on information derived from the previous solutions. This is achieved either by adapting the initial guess of the subsequent system or through a so-called subspace recycling, which allows constructing progressively more efficient two-level preconditioners. We report an overall speed-up over solving the systems independently of a factor of nearly $7$, or $5$, in our numerical experiments, which are inspired by the likelihood maximization and likelihood sampling procedures, respectively.
}


\maketitle

\section{Context and motivation}
\label{sect:intro:context}

Measurements registered by cosmic microwave background (CMB) experiments contain, in addition to the sought-after signal of cosmological origin, contributions from astrophysical sources. These are generically called foregrounds and can be of either galactic or extragalactic origins and be either diffuse or point-source-like morphologically.  A separation of the foreground signals from each other and, specifically, from the CMB signal is therefore an essential step of any modern CMB data analysis. This step is referred to as component separation. It is performed by capitalizing on either different electromagnetic frequency dependence and/or statistical properties of different signals~\citep[e.g.,][and references therein]{planckCompSep2016}.  In polarization the foreground signals tend to dominate the CMB signal over a broad range of angular scales and observational frequencies. The next generation of CMB observatories will therefore be only capable of delivering its science in full if high-precision statistically sound and reliable component separation techniques and their numerically efficient implementations are available. 

Component separation is a nonlinear operation. Based on data measured at multiple different frequency bands, it aims to simultaneously recover the frequency dependence of the foregrounds as well as their spatial morphology. It is commonly performed in a pixel domain and uses maps of the sky estimated for each frequency band and their statistical uncertainties as inputs. These objects are assumed to have been obtained in a preceding step of the data analysis that is called map-making. 

For concreteness, in this work we focus on the so-called parametric component separation approach~\citep[e.g.,][]{Brandt1994, Eriksen2008, Stompor2009}, where the frequency-scaling relations for each considered sky component are assumed to be given up to a limited number of unknown parameters, called foreground spectral parameters. However, the numerical techniques discussed here are more general and should be found useful also in other component separation methods. 

The component separation is typically performed in two steps. In the first step, the spectral parameters, or more generally, the mixing matrix elements, are estimated from the data, and in the second step, they are used to recover maps of sky components from the frequency maps.
This approach is conceptually simple and potentially very efficient computationally. The input frequency maps preserve essentially all the information present in a typically much larger initial raw data set, and their smaller sizes make them easier to store and operate on. 
For the next generation of CMB experiments, we expect to have as many as $n_\mathrm{t} \sim \mathcal{O}(10^{13}-10^{15})$ raw measurements, but only $\npix \sim \mathcal{O}(10^5-10^8)$ sky pixels.

The pixel-domain component separation approaches can ensure satisfactory performance but require a sufficiently precise statistical description of the frequency maps. This has to be derived from the raw measurements, which we refer to hereafter as time-domain data. In practice, this is often difficult because storage and computational cycles are limited. A general full covariance matrix of a single frequency map contains $\npix^2 \sim   \mathcal{O}(10^{10}-10^{16})$ elements, which would need to be stored in memory. Computing these elements costs at least $\mathcal{O}(\lambda)\,\mathcal{O}(n_\mathrm{t})$ floating point operations (flops). (Here $\lambda$ is a time-domain noise correlation length and can reach many thousands of samples.) The full computations therefore quickly become prohibitively expensive. This is the case even if the explicit inversion of the covariance matrix is replaced by some iterative procedure, which typically requires $\mathcal{O}(n_\mathrm{iter}\,\npix^2)$ flops, where the number of iterations, $n_\mathrm{iter}$ is usually on the order of $10^2$. Consequently, the best way forward in practice may be to invoke some approximations. This is often problematic as well, however, because a successful approximation needs to ensure sufficient accuracy to avoid introducing systematic biases in the estimated foreground parameters and later also in the component maps.

A general solution to the problem would be to avoid relying on the frequency maps at all and to perform all the calculation directly on the time-domain data. This would typically require memory on the order of $\mathcal{O}(n_\mathrm{t})$ and $\mathcal{O}(p\,n_\mathrm{iter}\,n_\mathrm{t}\,\ln\,\lambda)$ flops. The prefactor $p$ is on the order of unity for a typical map-making run, but in our case, it can vary widely between a few tens and many thousands. This highlights the challenge faced by the proposed approach. We note that while this is certainly very demanding, it is not necessarily prohibitive. Some of the best optimized existing map-making codes can already perform many hundreds of runs, for instance, as required in massive Monte Carlo simulations. The proposed approach may not only be more robust, but may be the only way forward if significant time-domain noise correlations are present, $\lambda \gg 1$. This is commonly the case in the CMB experiments, in particular, those operating from the ground. 

In this work, we explore some of the avenues that might render this approach tractable. We first identify the main computation-heavy step that unavoidably appears in any implementation of this technique.  We then investigate  how it might be accelerated by employing better and more advanced methods and their implementations.

The plan of this paper is as follows. In Section~\ref{sec:problem} we present the component separation problem and the numerical challenges it poses. In Section~\ref{sec:solutionproc} we describe the proposed solution and in Section~\ref{sec:numexp} the results of the numerical tests. Section~\ref{sec:conclusion} provides a brief summary and outlines prospects. Material that is more technical in nature or that is added for completeness is as usual deferred to the appendices.

\section{Problem description and setting}
\label{sec:problem}

\subsection{Preliminaries.}
Hereafter, we consider polarized signals and assume,  for simplicity and definiteness, that in every sky pixel the signal is characterized by two Stokes parameters, $Q$ and $U$. Extensions to include total intensity, $I$, are straightforward. Consequently, hereafter, every considered sky map consists of two maps corresponding to the two Stokes parameters. They are concatenated in a single map vector,
\begin{eqnarray}
    \sig{v} & = & \signal{v_q, v_u} \equiv \begin{bmatrix}
       v_q \\ v_u
    \end{bmatrix} , \qquad v_q, v_u \in \mathbb{R}^{\npix}.
\end{eqnarray}
Hereafter, partial brackets, $\lceil\dots\rfloor$, denote a vertical object. Examples of the sky maps as discussed in the following are single-frequency maps storing information about the sky signal as observed at a given frequency, or single-component maps containing information about a sky signal of some specific physical origin. 
We refer to the ordering defined above as Stokes-wise because a complete sky map of one Stokes parameter is followed up by another. In addition, we also consider a pixel-wise ordering, which for single maps reads
\begin{eqnarray}
    \sig{v} & \equiv & \signal{v_{q}(1), v_{u}(1), \ldots, v_{q}(\npix), v_{u}(\npix)},
\end{eqnarray}
where the $Q$ and $U$ parameters of a signal in one pixel are stored consecutively and are followed by those in another.

The goal of the component separation procedure is to estimate all assumed sky component signals given multiple frequency data. Therefore, we commonly deal with multiple maps of the same type, such as multiple single-frequency maps or multiple single-component maps. We concatenate them in a single multifrequency or multicomponent vector. For definiteness, in this work we fix the number of components to $n_\mathrm{comp}=3$ and consider three different sky components: CMB, dust, and synchrotron. A multicomponent vector, $\sigg{s}$, therefore contains information about the $Q$ and $U$ Stokes parameters of all three components. Such a vector can be ordered in multiple ways. Most commonly, we assume that it is ordered either in a component-wise way, when
\begin{eqnarray}
\label{eq:signal}
    \sigg{s} & \equiv & \signal{ \sig{s}_{\mathrm{cmb}}, \sig{s}_{\mathrm{dust}}, \sig{s}_{\mathrm{sync}} }\nonumber\\ 
    & = & \signal{ s_{\mathrm{cmb}, q}, s_{\mathrm{cmb}, u}, s_{\mathrm{dust}, q}, s_{\mathrm{dust}, u}, s_{\mathrm{sync}, q}, s_{\mathrm{sync}, u} }
     \in \mathbb{R}^{6\npix},
\end{eqnarray}
or in a pixel-wise way, where for each pixel all Stokes parameters follow consecutively for all considered components, that is, 
\begin{eqnarray}
    \sigg{s} & \equiv & \big\lceil s_{\mathrm{cmb}, q}(1), s_{\mathrm{cmb}, u}(1), \ldots, s_{\mathrm{sync}, q}(1), s_{\mathrm{sync}, u}(1), \ldots\\
    & & \ \ s_{\mathrm{cmb}, q}(\npix), s_{\mathrm{cmb}, u}(\npix),\ldots, s_{\mathrm{sync}, q}(\npix), s_{\mathrm{sync}, u}(\npix)\big\rfloor.
    \nonumber
\end{eqnarray}
Multifrequency vectors can be ordered in analogous manners.

The choice of the ordering in general depends on the specific context and is obviously of key importance for the numerical implementation of the map-making or component separation procedures.  Nonetheless, mathematically, switching the ordering from one to another is described by a linear, orthonormal, full-rank operator, $U$.  This operator is conceptually trivial to apply, and its application commutes with other matrix operations such as a matrix inversion because\begin{eqnarray}
(U\,M\,U^t)^{-1} & = & U\,M^{-1}\,U^t,
\end{eqnarray}
for any invertible matrix $M$. Consequently, a matrix can be inverted using one ordering, for instance, computing $M^{-1}$, and the result can later be reordered to obtain the inverse in the other ordering scheme, that is, $(U\,M\,U^t)^{-1}$. For this reason, we freely switch between the different orderings depending on the context in the following in order to highlight specific structures of the matrices, which may be more apparent for one choice than the other.

\subsection{Data model.}

As mentioned earlier, we consider a component separation procedure performed directly on the time-domain data as measured by the instrument. Thus we do not invoke any prior explicit map-making procedure. We therefore need to relate the time-domain measurements directly to the component maps because these maps are the intended outcome of the component separation procedure. We assume that for each frequency the time-domain data are made of sequences of consecutive observations registered by all detectors operating at this frequency and concatenated, we can write
\begin{align}
\label{eq:observed}
    d_f &= P_{\beta^\star, f} \, \sigg{s}_\star + n_f, \qquad d_f, n_f \in \mathbb{R}^\nt, \qquad f = 1, \ldots, \nfreq.
\end{align}
Here $\sigg{s}_\star$ is the unknown vector of the component amplitudes, and the star indicates that those are their actual values. $n_f$ denotes an (unknown) noise vector. The number of the frequency channels, $\nfreq$, is assumed to be larger than that of the components, $n_\mathrm{comp}$, set to $3$ in this work, to ensure that the problem is well defined. The matrix $P_{\beta^\star, f}$ in Eq.~(\ref{eq:observed}) combines the information about the instrument operations and the sky properties.  It can be expressed as
\begin{eqnarray}
    P_{\beta^\star, f} = P_f \cdot M_{\beta^\star, f},
    \label{eq:pDecomp}
\end{eqnarray}
where $M_{\beta^\star, f}  \in \mathbb{R}^{2\npix \times 6\npix}$ is  a so-called \emph{\textup{mixing matrix,}} and it determines how different sky components mix at all observed frequencies to yield the observed signal. The mixing matrix explicitly depends on the foreground scaling parameters, which we denote as $\beta^\star$, and the frequency of the observation,~$f$. $P_f \in \mathbb{R}^{\nt \times 2\npix}$ is in turn a \emph{\textup{pointing matrix}} defining which pixel of the sky each detector operating at a given frequency observed at every time. While it does not explicitly depend on frequency or scaling parameters, it therefore is in principle different for different frequencies because it encodes pointing of detectors specific to this frequency. This is highlighted by the subscript $f$. We have
\begin{align}
    P_f: \signal{s_{f,q}^\star, s_{f,u}^\star} & \mapsto d_f, \qquad M_{\beta^\star, f} : \sigg{s}_\star \mapsto \sig{s}_f \equiv \signal{s_{f,q}^\star, s_{f,u}^\star},
\end{align}
where $\sig{s}_f^\star $ is a single-frequency map expressing the combined sky signal at frequency $f$. The data vector, $d_f$, is time-ordered because its elements are indexed by the time at which the measurement was taken.

\subsection{Component separation.}

The goal of the component separation procedure is to solve an inverse problem, Eq.~(\ref{eq:observed}), and estimate the components, $\sigg{s}_\star$, given the full data set, $d \; (:= \, \{d_f\})$, made of data taken at all observational frequencies. This is typically solved by assuming that the noise, $n_f$, is Gaussian, with a zero mean and a known variance, $N_f$, and writing a data likelihood,
\begin{eqnarray}
-2\ln \mathcal{L}(\beta, \sigg{s}; d) & = & 
(\widetilde{d}-\widetilde{P}_{\beta} \,\sigg{s})^\top\,N^{-1}\,(\widetilde{d}-\widetilde{P}_{\beta} \,\sigg{s})\,+\,\hbox{\sc const}.
\label{eq:dataLike}
\end{eqnarray}
Here we have dropped $\text{the star}$ to distinguish an estimate from the true value, and we have introduced a tilde to denote multifrequency objects. We have
\begin{align}
    \widetilde{P}_{\beta} = \begin{bmatrix} P_{\beta, 1} \\ \vdots \\ P_{\beta, \nfreq}
    \end{bmatrix} = \begin{bmatrix} P_1 \cdot M_{\beta, 1} \\ \vdots \\ P_{\nfreq} \cdot M_{\beta, \nfreq}
    \end{bmatrix},
\end{align}
which follows from Eq.~(\ref{eq:pDecomp}), and
\begin{align}
    \widetilde{N}= \begin{bmatrix} N_1 & & 0  \\ & \ddots \\ 0 & & N_{\nfreq}
    \end{bmatrix},
    \quad
    \widetilde{d} = \begin{bmatrix}d_1\\ \vdots \\ d_{\nfreq}
    \end{bmatrix},
\end{align}
 which assumes no noise correlations between different frequency channels. In addition, throughout this work we also assume that while the component mixing represented by $M_\beta$ may involve (potentially) all components, it is always done on a pixel-by-pixel basis, so that all the elements of $M_\beta$ corresponding to different pixels vanish. Similarly, and in agreement with assumptions made in map-making procedures, we assume that the noise matrices, $N_f$, are block diagonal, with each block representing a banded Toeplitz matrix.

The standard two-step component separation procedure proceeds by first estimating for each frequency band, $f$, a single-frequency map, $m_f$,  and its covariance, $\hat{N}_f$. These are given by
\begin{eqnarray}
m_f & = & (P^\top_f\,N_f^{-1}\,P_f)^{-1}\,P^\top_f\,N_f^{-1}\,d_f, 
\label{eq:mapSolve}
\\
\hat{N}_f & = & (P^\top_f\,N_f^{-1}\,P_f)^{-1}.
\label{eq:covMap}
\end{eqnarray}
The follow-up component separation step is then performed assuming that the single-frequency maps yielded by the first step can be represented as
\begin{eqnarray}
m_f & = & M_{\beta^\star, f}\, \sigg{s}_\star \, + \, \hat{n}_f,
\label{eq:mapModel}
\end{eqnarray}
where $\hat{n}_f$ stands for a pixel-domain noise and is a Gaussian variable with variance $\hat{N}_f$. We can therefore write the corresponding likelihood as\begin{align}
-2\ln &\, \mathcal{L}(\beta, \sigg{s}; \{m_f\}) = \nonumber \\ & =  
\sum_f\,(m_f-M_{\beta, f} \,\sigg{s})^\top\,\hat{N}_f^{-1}\,(m_f-M_{\beta, f} \,\sigg{s})\,+\,\hbox{\sc const}.
\label{eq:dataLikePix}
\end{align}
This procedure is equivalent to directly solving the maximum likelihood problem defined by Eq.~\eqref{eq:dataLike}. However, it requires an explicit calculation of $\hat{N}_f^{-1}$ that for the current and forthcoming experiment is typically prohibitive because of restrictions on both the available computer memory and computational cycles. An alternative might be solving the original problem directly without explicitly invoking any pixel-domain objects. 
This is the option we study in this work. We note here in passing that intermediate approaches are also possible: for instance, one that relies on the likelihood in Eq.~\eqref{eq:dataLikePix}, but does not assume that $\hat{N}_f^{-1}$ is given explicitly. Instead, it computes a product of the covariance and a vector using an iterative procedure, which only requires applying the inverse covariance to a vector. This is performed using its implicit representation, Eq.~\eqref{eq:covMap}, as is done in the map-making solvers.
On the algorithmic level, such approaches are  equivalent to solving the problem in the time domain, and the methods considered hereafter would be applicable to that approach as well.

To estimate $\beta$ and $\sigg{s}$ directly from Eq.~\eqref{eq:dataLike}, we may either maximize this likelihood  or sample from a~posterior derived from it assuming some priors on the spectral parameters\footnote{We note that in sampling from the posterior, some priors for the signal would typically also be postulated, which would lead to a different system of equations than the one studied in this work. We leave this case to follow-up work.}. Alternatively,  a so-called spectral likelihood may be used~\citep[][]{Stompor2009}, where $\sigg{s}$ is already either marginalized or maximized over, that is,
\begin{align}
2\ln &\, \mathcal{L}_{spec}(\beta; \widetilde{d}) = \nonumber \\
&= \widetilde{d}^{\,\top}\,\widetilde{N}^{-1}\,\widetilde{P}_{\beta}\,(\widetilde{P}^\top_{\beta}\,\widetilde{N}^{-1}\,\widetilde{P}_{\beta})^{-1}\,\widetilde{P}^\top_{\beta}\,\widetilde{N}^{-1}\,\widetilde{d}\,+\,\hbox{\sc const},
\label{eq:specLike}
\end{align}
which again can be either minimized or sampled from. 

In both these cases, a key operation is a solution of a linear system of equations given by
\begin{equation}
\label{eq:origproblem}
    \widetilde{P}^\top_{\beta_i} \widetilde{N}^{-1} \widetilde{P}_{\beta_i} \, \sigg{s}_{\beta_i} = \widetilde{P}^\top_{\beta_i} \widetilde{N}^{-1} \widetilde{d},
\end{equation}
for a sequence of tentative values of the spectral parameters, $\beta_i$. These can be either a chain produced as a result of sampling, or a sequence of values obtained in the course of a minimization. We note that Eq.~(\ref{eq:origproblem}) is essentially a so-called map-making equation~\citep[e.g.,][]{Natoli2001, SzyGriSto14, Pugetal18}, but with a pointing matrix now replaced by $\widetilde{P}_{\beta_i}$. We can thus hope that as in the map-making problem, we can capitalize on special structures of the involved matrices and very efficient iterative solvers for solving linear systems to render the problem feasible. We point out that in the applications considered here, a subsequent value of the parameter $\beta$, that is, $\beta_{i+1}$, can only be known after the system for the current value, $\beta_i$, is fully resolved. A simultaneous computation of all systems for all values of $\beta_i$ is therefore not possible, and any computational speedup has to come from using better solvers for the linear systems and/or their numerical implementations.

When we separate parts that are dependent and independent of $\beta$, the system in Eq.~\eqref{eq:origproblem} can also be written as\begin{multline}
\label{eq:origproblem2}
\begin{bmatrix} 
    M_{\beta, 1} \\ \vdots \\ M_{\beta, n_{\mathrm{freq.}}}
\end{bmatrix}^\top
\overbrace{
\begin{bmatrix} 
    P_1^\top N_1^{-1} P_1 & & 0  \\ & \ddots \\ 0 & & P_\nfreq^\top N_\nfreq^{-1} P_\nfreq
\end{bmatrix}
}^{\let\scriptstyle\textstyle
\substack{\equiv \widetilde{A}}}
\overbrace{
\begin{bmatrix} 
    M_{\beta, 1} \\ \vdots \\ M_{\beta, n_{\mathrm{freq.}}}
\end{bmatrix}
}^{\let\scriptstyle\textstyle
\substack{\equiv \widetilde{M}_{\beta}}}
\sigg{s}_\beta
= \\ =
\begin{bmatrix} 
    M_{\beta, 1} \\ \vdots \\ M_{\beta, n_{\mathrm{freq.}}}
\end{bmatrix}^\top
\underbrace{
\begin{bmatrix} 
    P_1^\top N_1^{-1} d_1 \\ \vdots \\ P_\nfreq^\top N_\nfreq^{-1} d_{n_{\mathrm{freq.}}}
\end{bmatrix}
}_{\let\scriptstyle\textstyle
\substack{= \widetilde{P}^\top \widetilde{N}^{-1} \widetilde{d}}}
\end{multline}

The approach we propose here is based on two observations. First, our system has some essential similarities to that of the map-making problem, we should therefore be able to capitalize on novel iterative techniques proposed in that case. Second, we expect that consecutive values of $\beta_i$ in realistic sequences should not vary arbitrarily, and therefore
 subsequent linear systems~\eqref{eq:origproblem} should show some resemblance. Consequently, it should be possible to shorten the time to solution for the next value of $\beta_{i+1}$ by capitalizing on the solution for the current one, $\beta_i$.

\subsection{Block-diagonal preconditioner}

The block-diagonal preconditioner is the most common preconditioner used in the preconditioned conjugate gradient solvers applied in the context of the CMB map-making problem~\citep{Natoli2001}, which has demonstrated a very good performance in a number of applications. It is also the basis for the construction of more advanced preconditioners~\citep[e.g.,][]{SzyGriSto14}.
The block-diagonal preconditioner is derived by replacing the noise covariance~$N_f^{-1}$ in Eq.~\eqref{eq:covMap} by its diagonal. In the map-making case, when pixel-wise ordering is assumed, this leads to a block-diagonal matrix with the blocksize defined by the number of the considered Stokes parameters.
 In the component separation case, this preconditioner is given by $\widetilde{P}^{\,\top}_{\beta} \, \mbox{diag}(\widetilde{N}^{-1})\, \widetilde{P}_{\beta}$ , and in the pixel-wise ordering, it is block-diagonal. The diagonal block size is now equal to the product of the number of Stokes parameters and the number of sky components, that is, $ 6 \times 6$ in the specific case considered here. Consequently, the preconditioner can easily be inverted in any ordering scheme adapted.

Hereafter, we denote the $\beta$-independent part of the preconditioner as $\widetilde{B} := \widetilde{P}^\top \, \mbox{diag}(\widetilde{N}^{-1})\, \widetilde{P}$ so that
\begin{equation}
    \widetilde{P}^\top_{\beta} \, \mbox{diag}(\widetilde{N}^{-1})\, \widetilde{P}_{\beta} = \widetilde{M}_{\beta}^\top \, \widetilde{B}\, \widetilde{M}_{\beta}.
\end{equation}
By preconditioning the system~\eqref{eq:origproblem} from the left, we obtain
\begin{equation}
\label{eq:precproblem}
    \big( \widetilde{M}_{\beta}^\top \, \widetilde{B} \widetilde{M}_{\beta} \big)^{-1} \widetilde{M}_{\beta}^\top \, \widetilde{A} \widetilde{M}_{\beta} \, \sigg{s}_\beta
    = \big( \widetilde{M}_{\beta}^\top \, \widetilde{B} \widetilde{M}_{\beta} \big)^{-1} \widetilde{M}_{\beta}^\top \, \widetilde{P}^\top \widetilde{N} \widetilde{d}.
\end{equation}
To simplify the notation in the following, we define
\begin{eqnarray}
\A := \widetilde{M}_{\beta}^\top \, \widetilde{A}\, \widetilde{M}_{\beta}, \ \ \ \ \ \  \M := \widetilde{M}_{\beta}^\top \, \widetilde{B}\, \widetilde{M}_{\beta}, \ \ \ \ \ \ \sigg{b}:= \widetilde{M}_{\beta}^\top \, \widetilde{P}^\top \widetilde{N}\, \widetilde{d}.
\label{eq:matDefCompSep}
\end{eqnarray}

\subsection{Component mixing}

For concreteness, we assume throughout the paper the following component mixing scheme:
\begin{equation}
\label{eq:assumpmixing}
\begin{array}{rl}
    s_{f,q} &= \alpha_{f,1} \, s_{\mathrm{cmb}, q} + \alpha_{f,2}(\beta_d) \, s_{\mathrm{dust}, q} + \alpha_{f,3}(\beta_s) \, s_{\mathrm{sync}, q}\,, \\
    s_{f,u} &= \alpha_{f,1} \, s_{\mathrm{cmb}, u} + \alpha_{f,2}(\beta_d) \, s_{\mathrm{dust}, u} + \alpha_{f,3}(\beta_s) \, s_{\mathrm{sync}, u}\,,
\end{array}
\end{equation}
which follows the standard assumptions that there is no $Q$ and $U$ mixing, and that the scaling laws for the Stokes parameters $Q$ and $U$ of each components are the same. 
In the component-wise ordering, such mixing corresponds to the mixing matrix of the form ($I$ is the identity matrix, $2\times 2$ in this case)
\begin{equation}
    M_{\beta, f} = 
    \left[
    \begin{array}{c @{}c @{}c @{}c @{}c @{}c}
        \alpha_{f,1}\,I & 0 & \alpha_{f,2}(\beta_d)\,I & 0 & \alpha_{f,3}(\beta_s)\,I & 0\\
        0 & \alpha_{f,1}\,I & 0 & \alpha_{f,2}(\beta_d)\, I & 0 & \alpha_{f,3}(\beta_s)\,I
    \end{array}\right].
\end{equation}
The coefficients $\alpha_{f,i}$ encode the assumed scaling laws of the CMB, $i=1$, dust, $i=2$, and synchrotron, $i=3,$ where the last two depend on unknown scaling parameters, $\beta_d$ and $\beta_s$.
This matrix can be rewritten with the help of the Kronecker product as
\begin{equation}
\label{eq:kronmixing}
    M_{\beta, f} = 
    \left[
    \begin{array}{c @{}c @{}c @{}c @{}c @{}c}
        \alpha_{f,1} & 0 & \alpha_{f,2}(\beta_d) & 0 & \alpha_{f,3}(\beta_s) & 0\\
        0 & \alpha_{f,1} & 0 & \alpha_{f,2}(\beta_d) & 0 & \alpha_{f,3}(\beta_s)
    \end{array}\right] \otimes I.
\end{equation}
Hereafter, we drop the explicit dependence of the mixing coefficients on $\beta$, denoting them simply as $\alpha_{f,k}$.

\section{Solution procedure for the parametric component separation problem}
\label{sec:solutionproc}

A complete solution to the component separation problem has to successfully address two aspects. First, it needs to propose an efficient approach to solving the sequences of linear systems as in Eq.~\eqref{eq:precproblem}. Second, it has to combine it with an optimized procedure for the efficient determination of the new values of the parameters~$\beta$. This study addresses the former problem and focuses on the solution of a sequence of linear systems obtained for some sequences of the spectral parameters. In order to provide a fair comparison of various proposed techniques, we generate a sequence~$\{\beta_i\}$ beforehand and therefore, unlike in the actual applications, in our experiments, $\beta_{i+1}$ is in fact independent of the results of the preceding solution. %computed approximation~$\sigg{s}^{(final)}_{\beta_i}$. 
This ensures that the performance of all the considered solvers is evaluated on the identical sequences of linear systems.
\\
The overall solution scheme we adapt here is then as follows:
\begin{description}
\item[0)] Initialize $\beta_0$ and $\sigg{s}^{(0)}_{\beta_0}$ (typically $\sigg{s}^{(0)}_{\beta_0} := 0$), set $i:=0$.
\item[1)] Given $\beta_i$ and the initial guess $\sigg{s}^{(0)}_{\beta_i}$, solve the preconditioned problem, Eq.~\eqref{eq:precproblem}, deriving the current approximation~$\sigg{s}^{(final)}_{\beta_i}$.
\item[2a)] Determine the new parameters~$\beta_{i+1}$.
\item[2b)] Compute a new deflation space for the system associated with~$\beta_{i+1}$ using a recycling technique (see details below). This should not involve the value of~$\beta_{i+1}$ so that this step can be made in parallel with \textbf{2a)}.
\item[3)] Compute the initial guess $\sigg{s}^{(0)}_{\beta_{i+1}}$.
\item[4)] Set $i := i+1$ and go to \textbf{1)}.
\end{description}

\noindent
In the subsections below, we discuss steps \textbf{1)}, \textbf{2b)}, and \textbf{3)} in more detail.

\subsection{PCG with deflation and two-level preconditioners}
\label{sec:deflatedPCG}

Although the block-diagonal preconditioner has been shown to ensure good performance in the map-making experience, it has been pointed out that even better performance can often be achieved by employing so-called two-level preconditioners~\citep{SzyGriSto14,Pugetal18}.
Such preconditioners are built from the block-diagonal preconditioner, constituting the first level, and the second level is constructed from a limited number of vectors that are to be \emph{\textup{deflated}} (i.e., suppressed in the operator) in order to accelerate the convergence. These vectors are typically taken to be approximate eigenvectors of the system matrix corresponding to its smallest eigenvalues, which often hamper the convergence of PCG with the block-diagonal preconditioner.

We start from the case of \emph{\textup{deflation}} for the (unpreconditioned) conjugate gradient (CG) method.  
CG applied to a linear system $\A\sigg{s}=\sigg{b}$ with a given initial vector $\sigg{s}^{(0)}$ and an initial residual~$\sigg{r}^{(0)} := \sigg{b}-\A\sigg{s}^{(0)}$  builds implicitly orthogonal (residuals) and $\A$-orthogonal (search directions) bases of the Krylov subspace,
\begin{equation}
    \mathcal{K}_j(\A, \sigg{r}^{(0)}) = \mbox{span}\{\sigg{r}^{(0)}, \A\sigg{r}^{(0)}, \A^2 \sigg{r}^{(0)}, \ldots, \A^{j-1} \sigg{r}^{(0)}\},
\end{equation}
and the $j$th CG approximation~$\sigg{s}^{(j)} \in \sigg{s}^{(0)} + \mathcal{K}_j(\A, \sigg{r}^{(0)})$ is determined by the orthogonality condition on the $j$th residual~$\sigg{r}^{(j)} := \sigg{b}-\A\sigg{r}^{(j)}$,
\begin{equation}
    \sigg{r}^{(j)} \perp \mathcal{K}_j(\A, \sigg{r}^{(0)}).
\end{equation}

For a given set of deflation vectors, that is, the vectors to be suppressed, we denote by $\mathcal{U}$ the subspace spanned by these vectors. 
The deflation techniques replace the original operator $\A:\mathbb{R}^n \rightarrow \mathbb{R}^n$ by a deflated operator $\widehat{\A}:(\mathbb{R}^n \setminus \mathcal{U}) \rightarrow (\mathbb{R}^n \setminus \mathcal{U})$. The approximation is then sought over the augmented subspace (see, e.g.,~\cite{Gauetal13}),
\begin{equation}
    \sigg{s}^{(j)} \in \widehat{\sigg{s}}_0 + \mathcal{K}_j\big(\widehat{\A}, \widehat{\sigg{r}}^{(0)}\big) \cup \mathcal{U}
,\end{equation}
and the $j$th residual is required to be orthogonal to $\mathcal{K}_j(\widehat{\A}, \widehat{\sigg{r}}^{(0)}) \cup \mathcal{U}$. This effectively prevents the solver from exploring the subspace $\mathcal{U}$.

An extension of this for the PCG with the (first-level) preconditioner $\M$ is straightforward because we can use the PCG to implicitly build the Krylov subspace $\mathcal{K}_j(\M^{-1}\A, \M^{-1}\sigg{r}^{(0)})$. In the considered application, the preconditioner~$\M$ is the block-diagonal preconditioner. There are many variants of two-level preconditioners, and we summarize them briefly in~\Cref{sec:deflationvariants}. A more thorough survey can be found in~\cite{vuik:deflation09}, for example.

Each iteration of a deflated (P)CG, that is, with or without the first level, is more costly
than a single iteration of a standard (P)CG. The additional cost primarily depends on the number of deflated vectors, that is, the dimension of~$\mathcal{U}$, but also on the deflation variant. Building the subspace~$\mathcal{U}$ typically requires some preliminary computations, which can be as costly as solving the system~\citep[see, e.g.,][]{SzyGriSto14,Pugetal18}. Another approach, applicable to the cases when multiple systems need to be solved, is to construct the vectors "on the fly" during the solution of the systems themselves, thus hiding the additional cost. This is the approach we detail in the next sections.

\subsection{Subspace recycling}
\label{sec:recycle}

Several constructions of the deflation space have been adapted to solving a sequence of linear systems, for instance, those of  \cite{saad00:deflated-cg}, \cite{desturler06:recycleGMRES}, \cite{desturler06:recycleMINRES}, \cite{desturler17:recycleMINRES2}, and \cite{jolivet16:bgcrodr}. In this work, where the system matrix is symmetric positive definite (SPD), we follow~\cite{saad00:deflated-cg}. We build a subspace $\mathcal{Z} \subset \mathcal{K}_j(\widehat{\A}, \widehat{\sigg{r}}^{(0)})$ by storing some of the vectors computed during the previous run of (P)CG solver and determine the slowest eigenvectors of the operator~$\A$ restricted on the subspace $\mathcal{U}\cup\mathcal{Z}$. These are taken as the deflation vectors for the next solution. The resulting algorithm is given in \Cref{sec:fullalg}.

We can determine the subspace $\mathcal{Z}$  using either the residual or the search direction vectors forming (assuming the exact arithmetic) the orthogonal or an $\widehat{\A}$-orthogonal basis of $\mathcal{K}_j(\widehat{\A}, \widehat{\sigg{r}}^{(0)})$. Following~\cite{saad00:deflated-cg}, we choose here to use the search direction vectors. We retain the first $\dim_p$ search direction vectors, where $\dim_p$ defines the dimension of the so-called recycle subspace. We use the first vectors because the orthogonality among the computed vectors is gradually lost in CG; it is therefore better preserved in the initial iterations.

The techniques for approximating $k$~eigenvectors of the operator over a given subspace are well established. Among them, we note the Ritz and harmonic Ritz projections, which are described in detail in \Cref{sec:Ritz}. They lead to solving a (generalized) eigenvalue problem of small dimension,  in our case, of $\dim(\mathcal{U}) + \dim(\mathcal{Z})$. While~\cite{saad00:deflated-cg} suggested using the harmonic Ritz projection, we found the Ritz projection slightly more efficient in our numerical experiments, and we therefore include this in the full algorithm described in~\Cref{sec:fullalg}.
In another difference with~\cite{saad00:deflated-cg}, we assemble the (small) generalized eigenvalue problem matrices in the harmonic Ritz projection using the matrix-matrix products (see \Cref{alg:subrec} in \Cref{sec:Ritz}) instead of the optimized algorithm from~\citep[][Section~5.1]{saad00:deflated-cg}. This is because we expect that the additional computational cost in our application is negligible and we therefore opted for simplicity.  
 
There is no general recommendation for the choice of the number of deflation vectors, $k$, and the dimension of the recycling subspace, $\dim_p$. Higher~$k$ may result in an increase of the overall number of matrix-vector products (the system matrix has to be applied to $k$~deflation vectors before the deflated (P)CG is started for each system) and high $\dim_p$ may cause numerical instabilities in solving the eigenvalue problems that determine the new deflation vectors. On the other hand, the low values of $k$ and $\dim_p$ may not speed up the process sufficiently. We test this numerically in \Cref{sec:numexp}.

\subsection{Effect of the eigenvalue multiplicity.}

\label{sect:impEigMult}

One limiting factor to the efficiency of this approach, and more generally, to the performance of any two-level preconditioner with the deflation space estimated using standard iterative techniques such as Arnoldi or Lanczos iterations, comes from a higher multiplicity of eigenvalues, that is,  multiple eigenvectors with the same corresponding eigenvalue. This can arise either as a result of some symmetries in the scanning strategies in the case of the map-making systems of equations, or as similarities in the noise covariances of the different single-frequency maps in the case of the component separation problem as studied here; see \Cref{sec:multiEigen} for an example. Admittedly, such symmetries and/or similarities are not typically expected in the cases of real data analysis, but they can arise in the cases of simulated data, in particular if simplifying assumptions are adopted in order to speed up and/or simplify the simulation process.

To shed light on this problem, we consider an SPD matrix $A$ and assume that~$\lambda$ is an eigenvalue with multiplicity higher than one. This means that there exists a subspace $V$ with $\dim(V) > 1$ such that
\begin{equation}
    Av = \lambda v, \qquad \forall v \in V.
\end{equation}
Let $w$ be an arbitrary vector used to initiate the construction of a Krylov subspace and $w_V$ its projection onto $V$, that is,
\begin{equation}
    w = w_V + w', \quad w_V \in V,\ w' \perp V.
\end{equation}
%This decomposition is unique for any normal matrix~$A$ (in particular, an SPD matrix is normal). 
Then
\begin{equation}
    A^\ell w = A^\ell w_V + A^\ell w' = \lambda^\ell w_V + A^\ell w', \qquad A^\ell w' \perp V .
\end{equation}
Therefore, the $j$th Krylov subspace satisfies
\begin{align}
    \mathcal{K}_j(A, w) &= \mbox{span}\{w, Aw, A^2 w, \ldots, A^{j-1} w\} \\
    &= \mbox{span}\{ w_V \} \cup \mbox{span}\{w', Aw', A^2 w', \ldots, A^{j-1} w'\}, \nonumber  
\end{align}
and the intersection of $\mathcal{K}_j(A, w)$ and $V$ is at most a one-dimensional subspace spanned by $w_V$,
\begin{equation}
    \mathcal{K}_j(A, w) \cap V = \mbox{span}\{ w_V \}.
\end{equation}

Consequently, methods based on the Krylov subspace approximation, therefore including Lanczos and Arnoldi iterations (see \Cref{sec:arnoldilanczos} for more details), can recover one vector at most from the subspace spanned by multiple eigenvectors with the same eigenvalue. This may not be sufficient to allow for a construction of an efficient two-level preconditioner, however, in particular if the eigenvalue with many corresponding eigenvectors happens to be small: with the single precomputed vector we can only deflate a one-dimensional subspace of the entire multidimensional space as spanned by all these eigenvectors, and the remainder may continue hampering the convergence. % constructed via the combination of Lanczos/Arnoldi method.

This problem can be overcome by using a more advanced eigenvalue solver that can detect and handle the higher multiplicity of eigenvalues. An efficient implementation of such a solver is for instance provided by the ARPACK library~\citep{arpack}. In this case, the preconditioner may need to be precomputed with the help of such advanced routines, instead of constructing it on the fly as proposed here. If the presence of the eigenvalue multiplicity and the corresponding eigenvectors is known ahead of time, these vectors can be accomodated in the on-the-fly procedure proposed here. This is indeed the case we have encountered in one of the test cases discussed below.

We point out that the multiplicity of the eigenvalues is in principle advantageous for the standard (P)CG. In exact arithmetic, the effect of the whole subspace might be then eliminated at the cost of a single iteration. This fact is often used in the analysis of preconditioning methods based on preconditioners shifting some, possibly many, eigenvalues to the same value.

Last but not least, we emphasize that an eigenvalue multiplicity implies neither any indistinguishability of the corresponding eigenvectors nor a presence of degenerate modes in the solution, at least as long as the eigenvalue is not (numerically) zero. If the eigenvalue is not zero, the multiplicity only tells us that components of the right-hand side of the linear system that belong to the subspace spanned by the corresponding eigenvectors are merely weighted by the inverse system matrix in exactly the same way.

\subsection{Choice of the initial guess}
\label{sec:initguess}

The simplest and standard way to solve the sequence is to run the PCG method with the initial guess set to zero. However, some evidence exists showing that at least in the map-making case, this may not always be the best choice~\citep{PapGriSto18}, in particular in cases with high signal-to-noise ratios. In the case of a sequence of linear systems, all the systems involve the same initial data set with the same signal and noise content.  Even in data with a low signal-to-noise ratio, it may therefore be expected that adapting the initial guess following previous results may speed the process up in an interesting way. Consequently, we explore here two alternatives and show by numerical experiments that they are indeed much more efficient.

\subsubsection{Previous solution as the initial guess}
\label{sec:continuation}

A natural idea is to run the PCG for the (new) problem corresponding to~$\beta_{i+1}$ starting with the computed approximation~$\sigg{s}^{(final)}_{\beta_i}$,
\begin{equation}
\label{eq:continuation}
    \sigg{s}^{(0)}_{\beta_{i+1}} := \sigg{s}^{(final)}_{\beta_i}\,.
\end{equation}
This can be in particular efficient when the parameters $\beta_{i}$ and $\beta_{i+1}$ do not significantly differ and it is expected that so do $\sigg{s}_{\beta_i}$ and $\sigg{s}_{\beta_{i+1}}$. 
% In the numerical experiments (see below), we indeed observe that the number of iterations needed to reach the convergence can be significantly reduced if Eq.~\eqref{eq:continuation} is applied.

\subsubsection{Adapted previous solution as the new initial guess}
\label{sec:adaptation}

Eq.~\eqref{eq:continuation} can be further adapted by capitalizing on the particular structure of the mixing matrix. To start, we rewrite Eq.~\eqref{eq:origproblem} as
\begin{equation}
\label{eq:adaptation1}
    %\widetilde{M}_{\beta}^\top \, \widetilde{A} \widetilde{M}_{\beta} \, \sigg{s}_\beta
    %= \widetilde{M}_{\beta}^\top \, \widetilde{P}^\top \widetilde{N}^{-1} \widetilde{d} 
    %\ \  \Longleftrightarrow \ \ 
    \widetilde{M}_{\beta}^\top \big( \widetilde{A} \widetilde{M}_{\beta} \, \sigg{s}_\beta
    - \widetilde{P}^\top \widetilde{N}^{-1} \widetilde{d} \big)
    = 0.
\end{equation}
If the matrix~$\widetilde{M}_{\beta}$ were square (and nonsingular), then
\begin{equation}
    \widetilde{M}_{\beta} \, \sigg{s}_\beta = \widetilde{A}^{-1} \widetilde{P}^\top \widetilde{N}^{-1} \widetilde{d}
\end{equation}
would be the vector independent of~$\beta$. The solution~$\sigg{s}_\beta$ might then be interpreted as the coefficients with respect to\ the basis given by the columns of~$\widetilde{M}_{\beta}$. Therefore we would have
\begin{equation}
    \sigg{s}_{\bar{\beta}} = (\widetilde{M}_{\bar{\beta}})^{-1} \widetilde{M}_{\beta} \, \sigg{s}_\beta
\end{equation}
for arbitrary~$\bar{\beta}$ (for which~$\widetilde{M}_{\bar{\beta}}$ is nonsingular).

In our case, matrix~$\widetilde{M}_{\beta}$  is rectangular of size $ 2\,\nfreq\,\npix \times 6\, \npix $ and has full column rank. When the number of frequencies $\nfreq$ is not significantly higher than 3, we can generalize the above idea and use as the initial guess for the new system the vector
\begin{equation}
\label{eq:adaptedguess}
    \sigg{s}^{(0)}_{\beta_{i+1}} :=  (\widetilde{M}_{\beta_{i+1}})^\dagger \widetilde{M}_{\beta_{i}} \, \sigg{s}^{(final)}_{\beta_{i}},
\end{equation}
where $M^\dagger$ is the (Moore--Penrose) pseudo-inverse of~$M$,
\begin{equation}
    M^\dagger \equiv \big( M^\top M \big)^{-1} M^\top.  
\end{equation}
We recall our assumption that~$M$ is of full column rank. Clearly, for $\beta_{i+1} = \beta_{i}$,
%Clearly, this vector satisfies a natural requirement: for $\beta_{i+1} = \beta_{i}$, 
\begin{equation}
    (\widetilde{M}_{\beta_{i+1}})^\dagger \widetilde{M}_{\beta_{i}} = I \quad \mbox{holds, and therefore} \quad (\widetilde{M}_{\beta_{i+1}})^\dagger \widetilde{M}_{\beta_{i}} \, \sigg{s}_{\beta_{i}} = \sigg{s}_{\beta_{i}}.
\end{equation}

Finally, we note that the computation of the vector in Eq.~\eqref{eq:adaptedguess} is very cheap because of the Kronecker structure \eqref{eq:kronmixing} of the matrices~$\widetilde{M}_{\beta}$. Writing $\widetilde{M}_{\beta} = K_\beta \otimes I$, we obtain
\begin{equation}
    (\widetilde{M}_{\beta_{i+1}})^\dagger \widetilde{M}_{\beta_{i}} 
    = \left((K_{\beta_{i+1}}^\top K_{\beta_{i+1}})^{-1} K_{\beta_{i+1}}^\top K_{\beta_{i}} \right) \otimes I,
\end{equation}
in other words, only the matrices of size $2\,\nfreq \times 6$ need to be handled, and the cost of the proposed adaptation is nearly negligible.

\section{Numerical experiments}
\label{sec:numexp}

\subsection{Simulated data}
\label{sec:simdata}

For our numerical tests we use a simulated data set composed of time-ordered multifrequency observations with a correlated, '1/f', noise. The characteristics of this data set are as follows.

\subsubsection{Pointing matrix}

We adopt the simple scanning strategy used in~\cite{PapGriSto18}. The entire time-ordered data set is composed of $\mathcal{O}(10^8)$ measurements per frequency and divided into four consecutive subsets. The pointing is assumed to be the same for each frequency. The underlying sky is pixelized using the Healpix pixelization scheme~\citep{Healpix2005} with the resolution parameter $n_\textrm{side}$ set to $1024$. The scan consists of repetitive scanning of a rectangular patch made of $256$ pixel rows and columns. The scanning is either horizontal, that is, along the pixel rows, for the first and third subset, or vertical, that is, along the pixel columns for the second and fourth subset. During a single left-to-right, or bottom-up sweep, each sky pixel is sampled only once, and the direction of the polarizer, $\varphi_t$, is fixed for each of the four subsets and is equal, with respect to the sky, to $0$, $\pi/4$, $\pi/2$, and, $3\pi/4$.

The sky signal contribution to every measurement is modeled as
\begin{eqnarray}
d_c(t) & = & Q_c^\star( p(t)) \cos 2\varphi_t \, + \, U_c^\star(p(t))\, \sin2\varphi_t,
\end{eqnarray}
where $p(t)$ denotes the pixel observed at time $t$, we do not include the total intensity, and $Q_c$ and $U_c$ stand for $Q$ and $U$ Stokes parameters of the combined, CMB + foregrounds,  sky signal observed at frequency $\nu_c$.

\subsubsection{Sky maps}

\label{sect:sky_maps_def}

We assume $\text{six}$ frequency channels that approximately correspond to those accessible for a ground-based experiment. These are
\begin{eqnarray}
\nu_c & \in & \Big\{ 30, 40, 90, 150, 220, 270 \Big\}\; \mathrm{GHz}.
\end{eqnarray}
The sky signal is composed of emissions from three sources: CMB, dust, and synchrotron. The CMB signal is simulated using the current best-fit CMB model~\citep{planckParams2016}, while we use the so-called COMMANDER templates~\citep{planckCompSep2016} to model the dust and synchrotron signals that we scale to our reference frequency, $\nu_\mathrm{ref} = 150$GHz, using Planck's fiducial laws.

For the scaling laws we take a blackbody for the CMB component ($T_\mathrm{CMB}=2.7525$K), a power law for the synchrotron, and a modified blackbody for the dust, therefore
\begin{eqnarray}
\mathcal{S}^\textrm{sync}(\nu, \beta_s^\star) & = & \nu^{\beta_s^\star}\\
\mathcal{S}^\textrm{dust}(\nu, \beta_d^\star, T_d^\star) & = & \left(\frac{h\nu}{k T_d^\star}\right)^{\beta_d^\star}\,B(\nu, T_d^\star),
\end{eqnarray}
where the star distinguishes the true values of the parameters, 
\begin{eqnarray}
\beta^\star & \equiv & \big[ \beta_s^\star, \beta_d^\star, T_d^\star\big] \; = \; \big[ -3.1, 1.59, 19.6\,\mathrm{K}\big],
\label{eq:specParsTrueVals}
\end{eqnarray}
and $B(\nu, T)$ denotes a blackbody at temperature, $T$. The simulated maps are expressed in thermodynamic units and
 are given by
\begin{multline}
Q_p^\star(\nu) = %\mathcal{S}^\mathrm{cmb}(\nu)\,
Q_p^{\mathrm{cmb},\, \star} \, + \, 
\Gamma_\mathrm{RJ}(\nu)\,\bigg[\frac{\mathcal{S}^\mathrm{dust}(\nu, \beta_d^\star, T_d^\star)}{\mathcal{S}^\mathrm{dust}(\nu_\mathrm{ref}, \beta_d^\star, T_d^\star)}
\,Q_p^{\mathrm{dust},\, \star}(\nu_\mathrm{ref}) \, \\
+ \, \frac{\mathcal{S}^\mathrm{sync}(\nu, \beta_s^\star)}{\mathcal{S}^\mathrm{sync}(\nu_\mathrm{ref})}
\,Q_p^{\mathrm{sync},\, \star}(\nu_\mathrm{ref}, \beta_s^\star)\bigg]
\end{multline}
for each frequency $\nu = \nu_c$ and each observed sky pixel $p$. An analogous formula holds for the Stokes $U$  parameter. Here, $\Gamma_\mathrm{RJ}(\nu)$ stands for a conversion factor from Rayleigh-Jeans to thermodynamic units. This expression is consistent with Eq.~\eqref{eq:kronmixing} upon a suitable definition of the coefficients $\alpha$.

In our numerical experiments, we fix the dust temperature, $T_d$, to its true value and assume that only the spectral indices, $\beta\,=\, [\beta_s, \beta_d]$, are determined from the data. We assume that these are estimated by maximizing the spectral likelihood, Eq.~\eqref{eq:specLike}, using a truncated Newton maximization procedure. We use this approach to generate a single sequence of $\{\beta_i\}$, which, as explained in Sect.~\ref{sec:solutionproc}, we adopt consistently in all our runs. The sequence is made of $26$ values and is shown in Fig.~\ref{fig:betas}. In Appendix D we show for completeness the results of a similar test, but performed for a sequence of $\beta$ derived by sampling of the spectral likelihood. The main conclusions derived in these two examples are consistent.

\subsubsection{Noise}

We assume a correlated noise in the time domain with a spectrum given by
\begin{eqnarray}
    P(f) & = & \sigma^2_{rms} \, (1+\frac{f_\mathrm{knee}}{f}),
\end{eqnarray}
where $f$ is the time-domain frequency. The values of $f_\mathrm{knee}$ adopted here are different for different frequency channels and taken to be such that there are strong noise correlations within a single sweep. They span the range from $0.5$ up to $3$Hz from the lowest to the highest frequency channel, respectively. The noise is apodized at very low frequencies, so that the noise power is finite. $\sigma^2_\mathrm{rms}$ is taken to be about $30\mu$K per sample, reflecting the fact that each measurement effectively corresponds to the combined measurement of many modern detectors operating at the same frequency. This leads to sky maps with a noise $\sigma_\mathrm{rms}^{Q/U} \sim 2-3\mu$K per pixel.

\begin{figure*}[ht]
    \centering
    \includegraphics[width=0.8\textwidth]{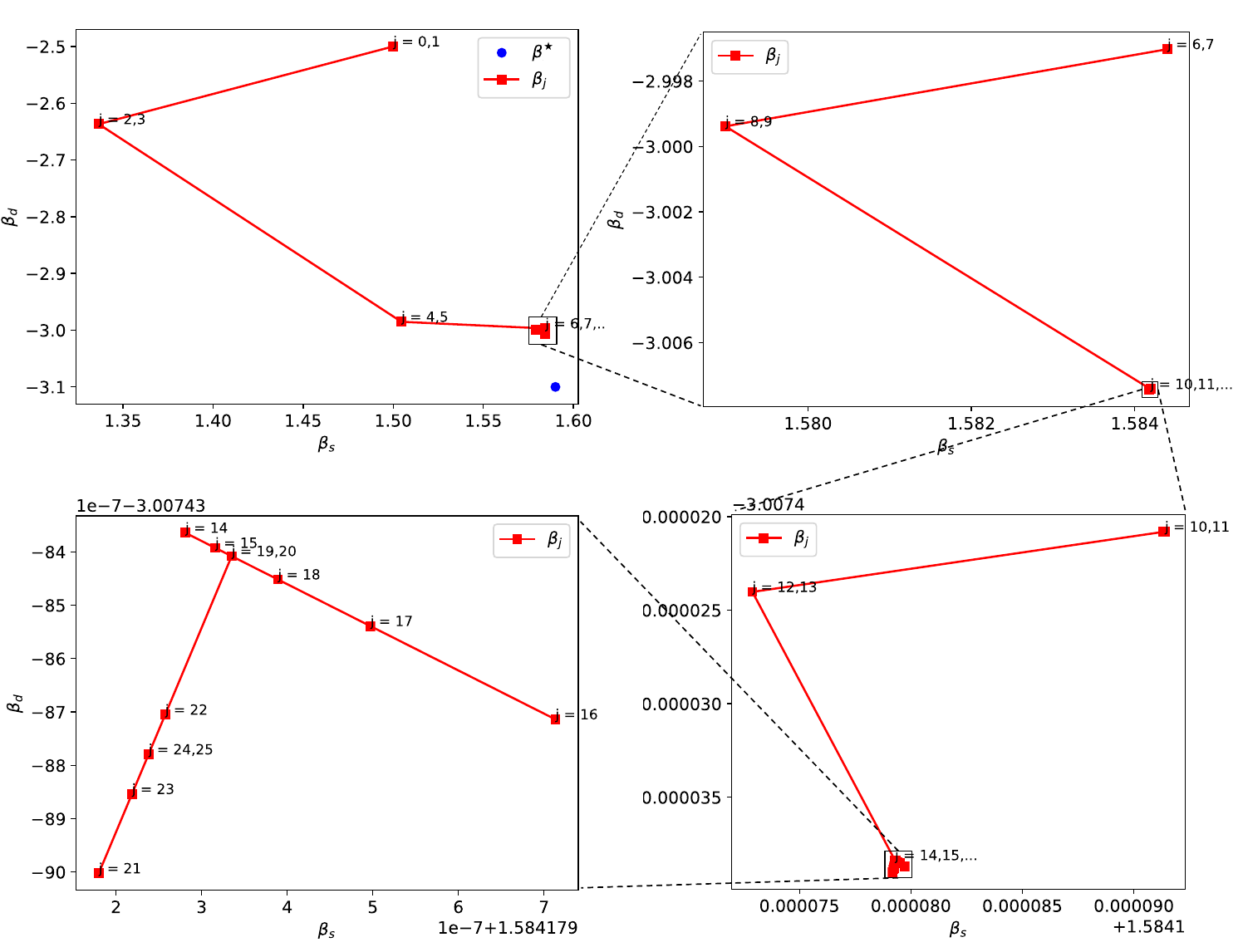}
    \caption{Sequence of the spectral parameters~$\beta_i = [\beta_{i,s}, \beta_{i,d}]$ used in our experiments and derived from the maximization of the spectral likelihood, Eq.~(\ref{eq:specLike}). The panels in clockwise order show consecutive zooms on the sequence that converged on the likelihood peak values of $\beta = [-3.006, 1.584],$ slightly off the true values marked by the blue solid point in the top left panel and given by Eq.~(\ref{eq:specParsTrueVals}) (with $T_d$ fixed in our test cases). The sequence was generated as described at the end of Sect.~\ref{sect:sky_maps_def}.}
    \label{fig:betas}
\end{figure*}

\subsection{Multiplicity of the eigenvalues as a result of the particular scanning strategy}
\label{sec:numexp_multieig}

In this section we address the issue of multiple eigenvectors with the same eigenvalues, which we have identified in our numerical tests. In agreement with Sect.~\ref{sect:impEigMult}, these were found to have significant effect on the performance of the proposed solvers. We show here how they can be traced back to the specific scanning strategy adopted in our simulations. We then describe corrective measures we included in order to minimize their effect and to ensure that our results are indeed representative of more typical cases. These considerations are given here for completeness and as a potentially useful example. However, as the proposed measures may not be necessary in most of the cases of realistic data, this section can be skipped without affecting the further reading.

We denote two pointing matrices for two horizontal scans as $P_0$ and $P_{\pi/2}$ and two pointing matrices for two vertical scans as $P_{\pi/4}$ and $P_{3\pi/4}$, where the subscript stands for the polarizer angle in the sky coordinates. They are related as 
\begin{equation}
\begin{array}{rcl}
    P_{\pi/2} & = & P_{0}\,
    \mathcal{R}_2(\pi/4),
    \\
 P_{3\pi/4} & = & P_{\pi/4}\,\mathcal{R}_2(\pi/4),  
\end{array}
\end{equation}
where $\mathcal{R}$ is a $12\,n_\mathrm{pix}$-by-$12\,n_\mathrm{pix}$ block-diagonal matrix with each diagonal block equal to a 2-by-2 spin-2 rotation matrix for each pixel of the $\text{six}$ frequency maps. While this is clearly a result of the simplifying assumption made in the simulations, this example may be of interest also in more realistic circumstances where certain relations of this sort may happen to be fulfilled approximately following from some common but inexact symmetries of typically adopted scans. We therefore briefly explore the consequences of this here.

In the case at hand, we can represent the total pointing matrix as
\begin{equation}
    \widetilde{P} =
    \left[
    \begin{array}{c}
         P_0\\
         P_{\pi/4}\\
         P_{\pi/2}\\
         P_{3\pi/4}
    \end{array}
    \right] \; = \; 
   \left[
    \begin{array}{c}
         P_{-\pi/4}\\
         P_{0}\\
         P_{\pi/4}\\
         P_{\pi/2}
    \end{array}
    \right] \, \mathcal{R}_{\pi/4} \, = \, 
    \left[
    \begin{array}{c}
         P_{3\pi/4}\\
         P_{0}\\
         P_{\pi/4}\\
         P_{\pi/2}
    \end{array}
    \right] \, \mathcal{R}_{\pi/4}
    \, = \, \widetilde{P}' \, \mathcal{R}_{\pi/4},
\end{equation}
given that all four scan subsets observe exactly the same sky.

When in addition the noise covariance for each of the four scan subsets is exactly the same, we have
\begin{equation}
    \widetilde{P}^\top\,\widetilde{N}^{-1}\,\widetilde{P} \, = \, \widetilde{P}'^\top\,\widetilde{N}^{-1}\,\widetilde{P}'
    \, = \,
    \mathcal{R}_{\pi/4}^\top\,\widetilde{P}^\top\,\widetilde{N}^{-1}\,\widetilde{P}\,\mathcal{R}_{\pi/4}.
\end{equation}
We note that this holds if the noise properties vary from one frequency channel to another, as is indeed the case in our simulations.

If now $v$ is an eigenvector of the matrix $\widetilde{A} = \widetilde{P}^\top\,\widetilde{N}^{-1}\,\widetilde{P}$ with a corresponding eigenvalue, $\lambda_v$, then
\begin{eqnarray}
  \widetilde{A}\,v & = & 
  \widetilde{P}^\top\,\widetilde{N}^{-1}\,\widetilde{P}\;v \; = \; \mathcal{R}_{\pi/4}^\top\,\widetilde{P}^\top\,\widetilde{N}^{-1}\,\widetilde{P}\;\mathcal{R}_{\pi/4} \, v \; = \; \lambda_v\, v,
  \label{eq:pmatSimsRel}
\end{eqnarray}
and therefore,
\begin{equation}
\widetilde{A} \;\mathcal{R}_{\pi/4} \,v \; = \;
\widetilde{P}^\top\,\widetilde{N}^{-1}\,\widetilde{P}\;\mathcal{R}_{\pi/4} \, v \; = \; \lambda_v\; \mathcal{R}_{\pi/4}\,v.
\end{equation}
This means that also $\mathcal{R}_{\pi/4}\,v$ is an eigenvector of the matrix $A$ with the same eigenvalue, $\lambda_v$. Because this reasoning applies as much to the matrix $A$ as the matrix $B = \widetilde{P}^\top\,\mathrm{diag}\,\widetilde{N}^{-1}\,\widetilde{P}$, we have
\begin{equation}
\begin{array}{rcl}
\widetilde{P}^\top\,\widetilde{N}^{-1}\,\widetilde{P}\;v & = & \lambda'_v\,\widetilde{P}^\top\,\mathrm{diag}\,\widetilde{N}^{-1}\,\widetilde{P}\;v\\
\widetilde{P}^\top\,\widetilde{N}^{-1}\,\widetilde{P}\;( \mathcal{R}_{\pi/4}\,v ) & = & \lambda'_v\,\widetilde{P}^\top\,\mathrm{diag}\,\widetilde{N}^{-1}\,\widetilde{P}\; (\mathcal{R}_{\pi/4}\,v ).
\end{array}
\end{equation}
In general, this does not yet imply that the component separation system matrix preconditioned with the block-diagonal preconditioner, Eq.~\eqref{eq:matDefCompSep}, given by
\begin{equation}
    \M^{-1}\,\A \; = \; (\widetilde{M}_\beta^\top\,\widetilde{B}\,\widetilde{M}_\beta)^{-1}
(\widetilde{M}_\beta^\top\,\widetilde{A}\,\widetilde{M}_\beta),
\end{equation}
has two eigenvectors corresponding to the same eigenvalue related by the rotation operator acting in the component space. 
This is the case when the subspace spanned by $v$ and $\mathcal{R}_{\pi/4}\,v$ is contained in the subspace spanned by the columns of the mixing matrix, $\widetilde{M}_\beta$. Otherwise, the preconditioned system matrix may have a single (when these two subspaces merely intersect) or no corresponding eigenvectors (when these two subspaces are disjoint). Which of these situations is actually realized is case dependent and in general also depends on the value of $\beta$.

We found that in our numerical experiments the eigenvalue multiplicity of the preconditioned system matrix due to the assumed scan symmetries was sufficiently common that we opted to account for it explicitly in our analysis. Consequently, we use the subspace recycling to approximate one of the eigenvectors, and we compute the other by applying the rotation operator. We then use both vectors to construct the deflation operator. Given that $\mathcal{R}_{\pi/4} = -\mathcal{R}_{-\pi/4}$, there is no ambiguity because we do not know a priori which of the two vectors we estimate directly through the subspace recycling, and this approach leads to the same deflation space, regardless of the rotation that is applied. This technique has led to a significant speed-up in the cases studied below.

\subsection{Results}
\label{sec:numexp_results}

We compare the convergence using the relative norm of the $j$th residual,
\begin{equation}
    \frac{\|\widetilde{M}_{\beta}^\top \, \widetilde{P}^\top \widetilde{N} \widetilde{d} - \widetilde{M}_{\beta}^\top \, \widetilde{A} \widetilde{M}_{\beta} \, \sigg{s}^{(j)}_\beta \|}{\| \widetilde{M}_{\beta}^\top \, \widetilde{P}^\top \widetilde{N} \widetilde{d} \|}.
\end{equation}
The iterations for each system are stopped when this value drops below tolerance $\mbox{TOL} = 10^{-8}$, but we always perform at least one iteration. 

\begin{figure}[ht]
    \centering
    \includegraphics[width=0.7\linewidth]{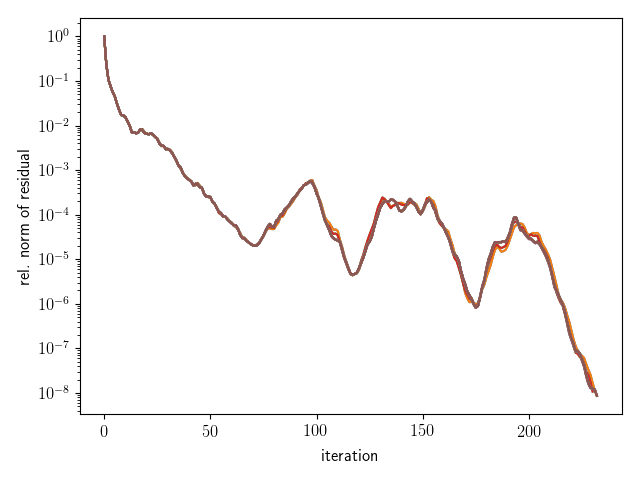}
    \caption{Convergence of PCG with zero initial guess and block-Jacobi preconditioner for all $26$ systems in the sequence shown in Fig.~\ref{fig:betas}.}
    \label{fig:allsystems_naive}
\end{figure}

We first show that the systems corresponding to different $\beta$s from the sequence are indeed "close to each other" in the sense that they display a similar behavior during the solution process. We illustrate this by showing the convergence of PCG with zero initial guess in \Cref{fig:allsystems_naive}. We find that all the convergence curves are indeed very similar, even for the initial elements of the sequence where the values of $\beta$ parameters continue to change quite significantly.

\begin{figure}[ht]
    \centering
    \includegraphics[width=0.7\linewidth]{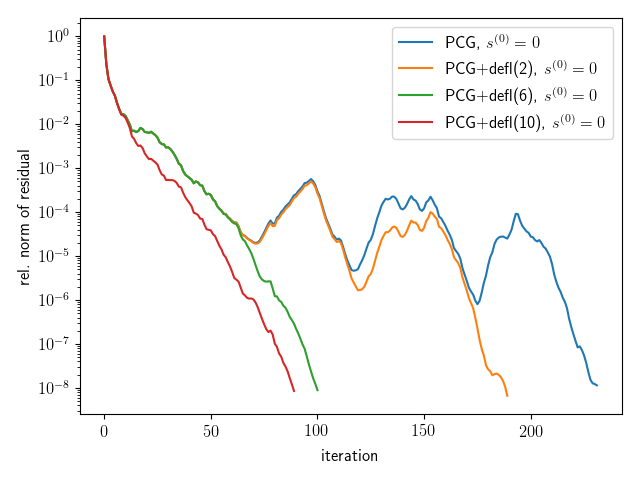}
    \caption{Convergence of PCG with the deflation applied to $2$, $6$, and $10$ slowest eigenvectors for the true system with $\beta = \beta^\star$.}
    \label{fig:defl_beta-1}
\end{figure}

We can therefore focus on the "true" system corresponding to $\beta = \beta^\star$ in order to investigate the improvement caused by deflating the eigenvectors corresponding to the smallest eigenvalues of the system matrix. This is depicted in \Cref{fig:defl_beta-1}. Here, the eigenvectors are computed using the ARPACK eigensolver \citep{arpack}, and as expected, we find that significant improvement is achieved by the deflation.

\begin{figure}[ht]
    \centering
    \includegraphics[width=0.7\linewidth]{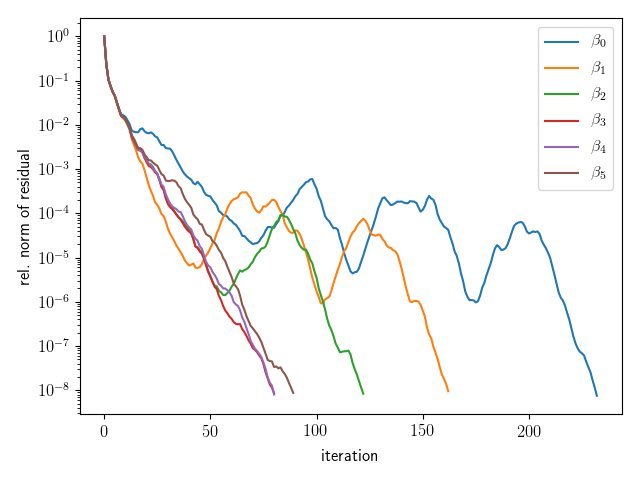}
    \caption{Convergence of deflated PCG with the deflation subspace build by recycling. Here we consider the first six systems from the sequence and start the iterations always with zero initial guess. $k=10$ eigenvectors are approximated using $\dim_p = 100$.}
    \label{fig:recy_only}
\end{figure}

Then, we illustrate the effect of the deflation space built by recycling. To do this, we first consider six systems and start the iterations always with zero initial guess, $\sigg{s}^{(0)} = 0$. The result for $k=10$ eigenvectors approximated using the dimension of the subspace, $\dim_p = 100$, is given in \Cref{fig:recy_only}.

\begin{figure*}[htp]
    \centering
    \includegraphics[width=0.9\textwidth]{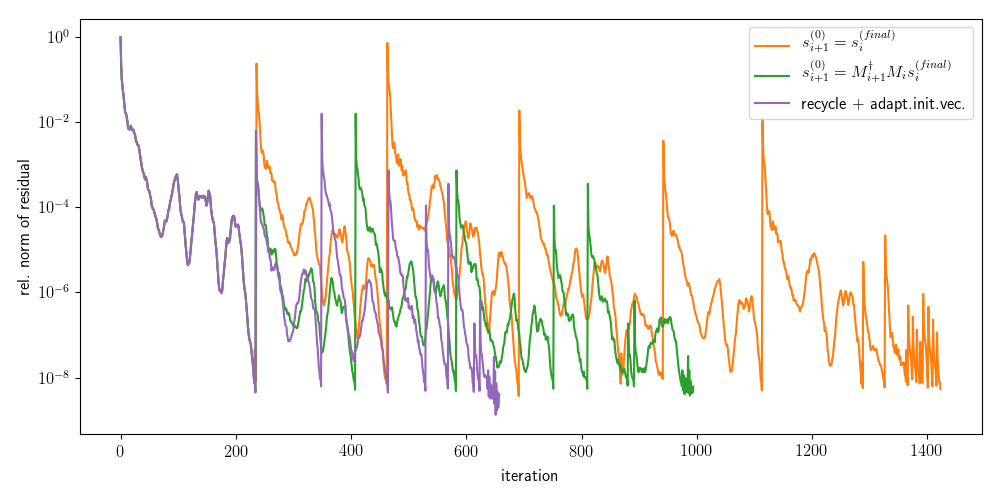}
    \caption{Comparison of the PCG with different choices of initial guess (as in Eq.~\eqref{eq:continuation} and~Eq.~\eqref{eq:adaptedguess}) and the PCG with the subspace recycling (together with the choice of the initial guess as in~Eq.~\eqref{eq:adaptedguess}). For the recycling we consider $k=10$ eigenvectors approximated using $\dim_p = 100$.}
    \label{fig:allsystems_compare}
\end{figure*}

In \Cref{fig:allsystems_compare} we compare the convergence of the full sequence of the $26$ systems using the techniques of setting the initial guess proposed in Eq.~\eqref{eq:continuation} and Eq.~\eqref{eq:adaptedguess} and using subspace recycling. We recall that standard PCG with zero initial vector takes more than $6000$ iterations; cf.\ Figure~\ref{fig:allsystems_naive}. Although the subspace recycling variant requires  $25\times 10$ matrix-vector products\footnote{It is necessary to apply the matrix to deflation vectors at the beginning of the deflated PCG to build the projection matrices, see \Cref{alg:deflPCG} in \Cref{sec:fullalg}.} more than the PCG with block-Jacobi preconditioner for any choice of the initial guess, it still provides an interesting increase in speed.

We also compare the time required by one iteration of the PCG with a deflation with the time required by a standard PCG iteration. In Table~\ref{tab:timings} we list the relative times of a single iteration of a deflated PCG with 2, 6, and 10 deflation vectors in our experiments (taking an average from five~runs, each with ten~iterations). The small and negligible overhead introduced by the two-level preconditioner indicates that most of the time in the code is spent on the standard map-making operations, such as (de)projection operators and noise weighting~(e.g., \citet{CanBorrJaf2010}). We emphasize that these timings depend on the implementation and the hardware on which it is executed. For massively parallel codes, the cost of a deflated PCG iteration compared to a standard PCG iteration may increase somewhat, but we expect that the proposed algorithm should nevertheless still be offering an attractive increase in speed.

\begin{table}[htp]
    \centering
    \begin{tabular}{cccc}
         PCG & PCG+defl(2) & PCG+defl(6) & PCG+defl(10)  \\ \hline
         1 & 1.0001 & 1.0013 & 1.0019\\
         ~
    \end{tabular}
    \caption{Timings of a single deflated PCG iteration with 2, 6, and 10 deflation vectors. The timings are relative with respect to a single iteration of the standard (nondeflated) PCG. The table gives the average from five~runs, each with ten~iterations.}
    \label{tab:timings}
\end{table}

\begin{figure*}[htp]
    \centering
    \includegraphics[width=0.9\textwidth]{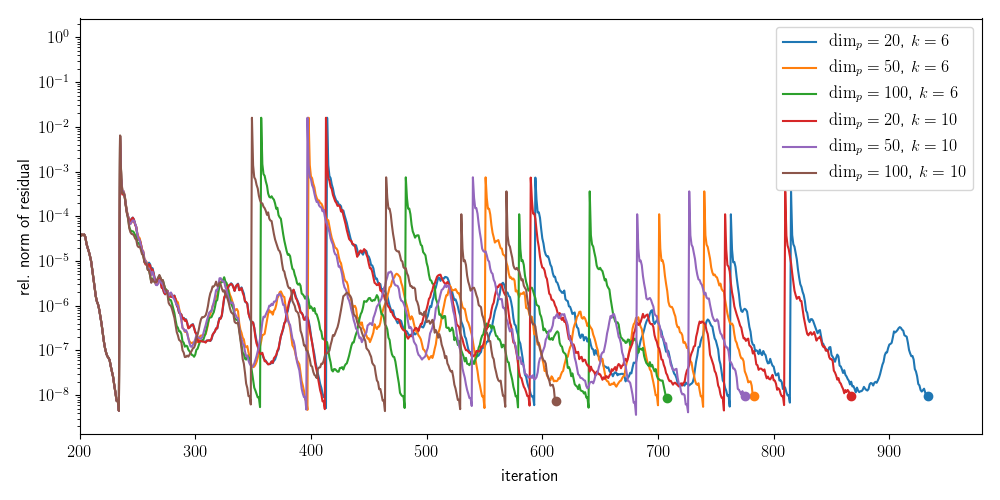}
    \caption{Comparison of the PCG with different choices of $k$ and $\dim_p$ for the first ten systems in the sequence. The initial guess is the same as in Eq.~\eqref{eq:adaptedguess}. The iteration counts are also listed in \Cref{tab:recy_var_k_dimP}. Because the convergence for the first system is independent of $\dim_p$ and $k$, the $x$-axis is depicted starting from the iteration 200.}
    \label{fig:recy_var_k_dimP}
\end{figure*}

As noted above, the performance of the PCG with the deflation space built by recycling is affected by the number of the deflation vectors, $k$, and the dimension of the recycling subspace, $\dim_p$. There is no general recommendation for their choice. Higher~$k$ may result in an increase of the overall number of matrix-vector products (the system matrix has to be applied to $k$~deflation vectors before the deflated PCG is started for each system) and high $\dim_p$ may cause numerical instabilities in solving the eigenvalue problems that determine the new deflation vectors. On the other hand, the low values of $k$ and $\dim_p$ may not increase the speed sufficiently. We compare the convergence of PCG with some choices of $k$ and $\dim_p$ in \Cref{fig:recy_var_k_dimP} and in \Cref{tab:recy_var_k_dimP}. The deflation clearly has a small effect for small $\dim_p$, that is, when the eigenvectors are not approximated accurately.

\begin{table}[htp]
\centering
\begin{tabular}{ccccc}
 $\dim_p$ & $k$ & \#iterations & \multicolumn{2}{c}{\#MatVecs}  \\ 
 & & & deflation & total \\ \hline
 20 & 6 & 933 & 54 & 987\\
 50 & 6 & 783 & 54 & 837\\
100 & 6 & 708 & 54 & 762\\
 20 & 10 & 867 & 90 & 957\\
 50 & 10 & 775 & 90 & 865\\
100 & 10 & 612 & 90 & 702
\end{tabular}
\caption{Number of PCG iterations and matrix-vector products (MatVecs) for different choices of $k$ and $\dim_p$ for the first $ten$ systems in the sequence. The initial guess is the same as in Eq.~\eqref{eq:adaptedguess}.}
    \label{tab:recy_var_k_dimP}
\end{table}

\section{Conclusions and further perspectives}
\label{sec:conclusion}

We have presented a procedure for efficiently solving a sequence of linear systems arising from the CMB parametric component separation problem. Two main novelties are the proposed choice of the initial vectors and the recycling technique used to determine the deflation space. 
Motivated by our simulated data, we also emphasized and elaborated on the role of the multiplicity of the eigenvalues, in particular in the context of their effect on the performance of two-level preconditioners.

The overall speed-up factor we obtained, $\sim$~5--7, is significant.  The bulk of the improvement comes from reusing the solution of an earlier system as the initial guess of the next solution - a simple but not trivial observation owing to the fact that this is the same data set being used in all the solutions. However, other proposed amendments provide a significant additional performance boost on the order of $\sim$~2. This is particularly significant in the case of the sampling-based application.
Further improvements and optimizations are clearly possible. For instance, the number of required matrix-vector products can be decreased by not using the two-level preconditioner for all the systems. As the experiments showed, when two consecutive system solutions are very similar, the PCG with the diagonal preconditioner and a proper initial guess (e.g.,~as proposed in \Cref{sec:adaptation}) can already be sufficient for convergence in a few iterations.

We emphasize that in practice, we will be only able to capitalize on this type of approach when they are implemented in a form of highly efficient high-performance massively parallel numerical algorithms and codes. We leave a full demonstration of this to future work, noting here only that many of the relevant techniques have been studied in the past and recent literature, showing that this should indeed be feasible~\cite[e.g.,][]{CanBorrJaf2010, SudBorrCant2011,SzyGriSto14,PapGriSto18,SelBaerEri2019}. 

The techniques we presented, rendered in the form of efficient high-performance codes, should allow for the efficient maximization of the data likelihood or the posterior distributions in the component separation problems and produce reliable sky component estimates for at least some of the forthcoming data sets. However, in the cases of sampling algorithms, when many thousand linear systems may need to be solved, this still remains to be demonstrated, and further improvements will likely be required. They will depend in general on specific details of the employed sampling technique, however, and we leave them here as future work.

The techniques discussed here can also be used in other problems in CMB data analysis that require solving a sequence of linear systems.
In particular, they should be directly relevant for applications that estimate instrumental parameters, which commonly have to be included in more realistic data models and estimated from the data.

The codes used in this work are available from a GitHub repository\footnote{https://github.com/B3Dcmb/Accelerated-PCS-solvers}.

\begin{acknowledgement}
We thank Olivier Tissot for insightful discussions and Dominic Beck and Josquin Errard for their help with numerical experiments.
The first two authors' work was supported by the NLAFET project as part of European Union's Horizon 2020 research and innovation program under grant 671633.
This work was also supported in part by the French National Research Agency (ANR) contract ANR-17-C23-0002-01 (project B3DCMB).
This research used resources of the National Energy Research Scientific Computing Center (NERSC), a DOE Office of Science User Facility supported by the Office of Science of the U.S. Department of Energy under Contract No. DE-AC02-05CH11231.
\end{acknowledgement}

\appendix

\section{Eigenvalue multiplicity in the component separation problem. A worked example.}
\label{sec:multiEigen}

In this section we discuss the eigenvector and eigenvalue structure of the preconditioned matrix defined in Eq.~\eqref{eq:precproblem} in the context of eigenvalue multiplicity in some specific setup that in particular assumes that the pointing matrix is the same for each frequency and that the noise has the same covariance (up to a scaling factor) for each frequency, that is,
\begin{equation}
    P_f = P,\qquad  N_f = N, \qquad f = 1, \ldots, n_{\mathrm{freq.}}.
\end{equation}

While these requirements are not very likely to be strictly realized in any actual data set, they can be fulfilled approximately, leading to near multiplicities of the eigenvalues. If these are not accounted for, they may be as harmful to the action of the preconditioner as the exact multiplicities.
Moreover, this worked example demonstrates that the component separation problem is in general expected to be more affected by this type of effect than the standard map-making solver, for instance, therefore emphasizing that due diligence is necessary in this former application.

First, let $(\lambda_i, \sig{v}_i) = (\lambda_i, \signal{v_{i,q}, v_{i,u}} )$ be an eigenpair of the map-making matrix, that is, there holds
\begin{equation}
\label{eq:GEPmapmaking}
    P^\top N^{-1} P \,\sig{v}_i = \lambda_i P^\top \, \mbox{diag}(N^{-1}) \,P \,\sig{v}_i.
\end{equation}
We note that
\begin{equation}
\label{eq:Mu}
    \widetilde{M}_{\beta} \signal{\sig{v}_i, \sig{0}, \sig{0}} = 
    \widetilde{M}_{\beta} \begin{bmatrix} v_{i,q} \\ v_{i,u} \\ 0 \\ \vdots \\ 0 \end{bmatrix} =
\begin{bmatrix} 
   \alpha_{1,1}\, v_{i,q} \\  \alpha_{1,1}\, v_{i,u} \\ \vdots \\ \alpha_{\nfreq,1}\, v_{i,q} \\ \alpha_{\nfreq,1}\, v_{i,u}
\end{bmatrix}
=
\begin{bmatrix} 
   \alpha_{1,1}\, \sig{v}_i \\ \vdots \\ \alpha_{\nfreq,1}\, \sig{v}_i
\end{bmatrix}
\end{equation}
because of the form of the mixing we assumed in Eq.~\eqref{eq:assumpmixing}. 
Consequently, using Eq.~\eqref{eq:GEPmapmaking},
\begin{multline}
    \widetilde{A} \widetilde{M}_{\beta} \signal{\sig{v}_i, \sig{0}, \sig{0}} 
    =
    \begin{bmatrix} 
       \alpha_{1,1}\, P^\top N^{-1} P \,\sig{v}_i \\ \vdots \\ \alpha_{\nfreq,1}\, P^\top N^{-1} P \,\sig{v}_i
    \end{bmatrix}
    \\
    = 
    \begin{bmatrix} 
       \alpha_{1,1}\, \lambda_i P^\top \, \mbox{diag}(N^{-1}) P \, \sig{v}_i \\ \vdots \\ \alpha_{\nfreq,1}\lambda_i P^\top \, \mbox{diag}(N^{-1}) P \,\sig{v}_i
    \end{bmatrix}
    = \lambda_i \widetilde{B} \widetilde{M}_{\beta} \signal{\sig{v}_i, \sig{0}, \sig{0}}.
\end{multline}
Because the matrix $ \widetilde{M}_{\beta}^\top \widetilde{A} \widetilde{M}_{\beta} $ is assumed to be nonsingular (equivalently, because $\widetilde{M}_{\beta}$ is of full column rank), we can multiply this equation from the left by $\widetilde{M}_{\beta}^\top$ showing that $(\lambda_i,\signal{\sig{v}_i, \sig{0}, \sig{0}} )$ is the eigenpair of the matrix $ ( \widetilde{M}^\top_{\beta} \widetilde{B} \widetilde{M}_{\beta} )^{-1} \widetilde{M}_{\beta}^\top \widetilde{A} \widetilde{M}_{\beta} $.
We can proceed analogously for the vectors $\signal{\sig{0}, \sig{v}_i, \sig{0}}$ and $\signal{\sig{0}, \sig{0}, \sig{v}_i}$, with replacing in Eq.~\eqref{eq:Mu} $\alpha_{f,1}$ by~$\alpha_{f,2}$ and~$\alpha_{f,3}$, respectively.

There are $2 \, \npix$ eigenpairs for $(P^\top \, \mbox{diag}(N^{-1}) P)^{-1} (P^\top N^{-1} P)$. As we showed above, each of them generates three eigenpairs (with the same eigenvalue) of $  ( \widetilde{M}^t_{\beta} \widetilde{B} \widetilde{M}_{\beta} )^{-1} \widetilde{M}_{\beta}^\top \widetilde{A} \widetilde{M}_{\beta} $. This gives together $6 \, \npix$ eigenpairs, in other words, we have described the full spectrum of the preconditioned system matrix in Eq.~\eqref{eq:precproblem}.

Finally, we remark that all the eigenpairs of the preconditioned matrix  in the simplified setting are independent of the parameters~$\beta$. In this case, we suggest using a specialized eigensolver (e.g., ARPACK, \cite{arpack}) to compute the eigenpairs from Eq. \eqref{eq:GEPmapmaking}, build the triplets of eigenvectors $\signal{\sig{v}_i, \sig{0}, \sig{0}}$, $\signal{\sig{0}, \sig{v}_i, \sig{0}}$,  and $\signal{\sig{0}, \sig{0}, \sig{v}_i}$, and then use the deflated PCG with these vectors.

\begin{figure*}[ht]
    \centering
    \includegraphics[width=0.9\textwidth]{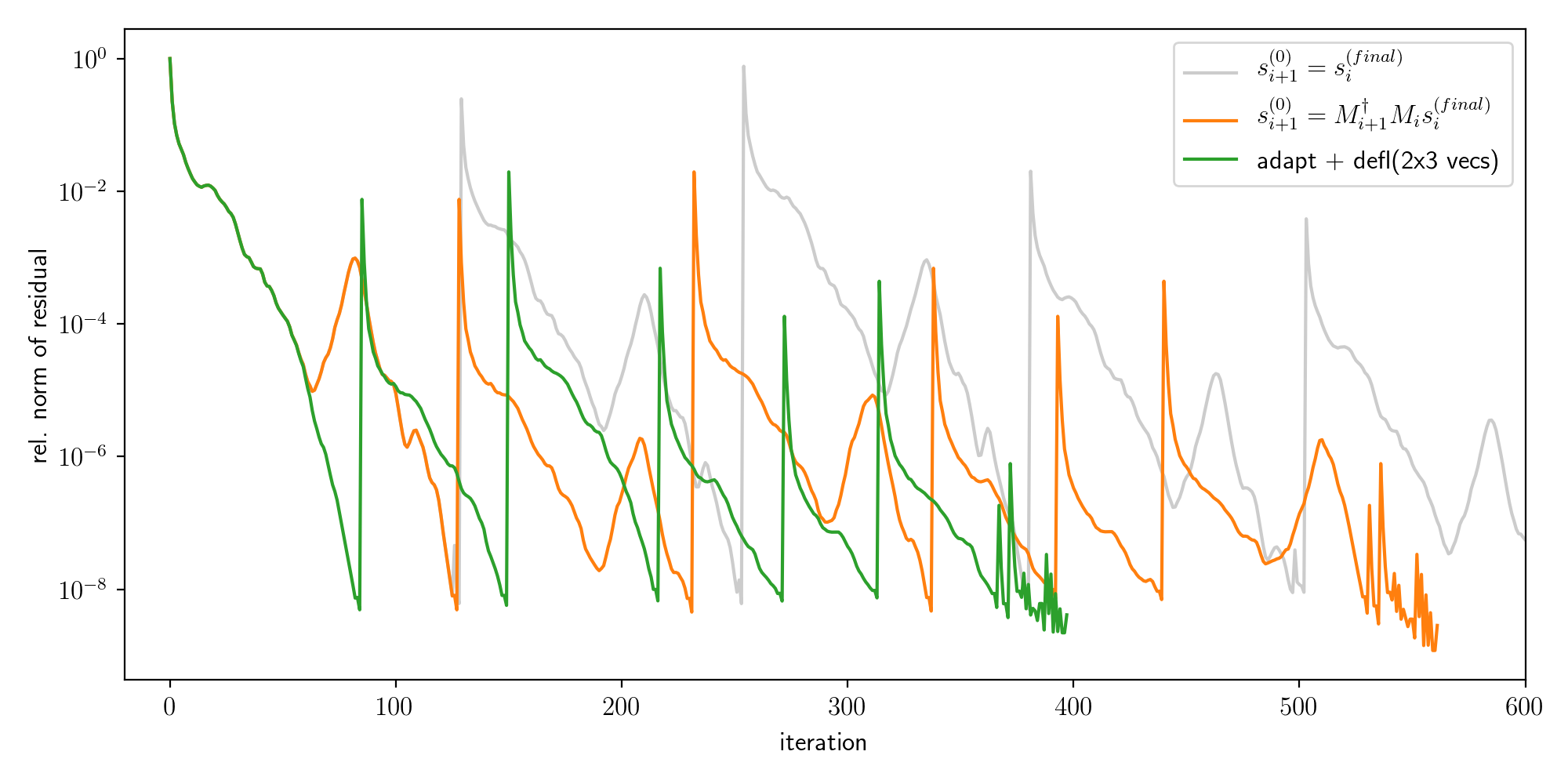}
    \caption{Same as \Cref{fig:allsystems_compare} for the simplified setting of Appendix~\ref{sec:multiEigen}. Comparison of the PCG with different choices of the initial guess (as in Eq.~\eqref{eq:continuation} and Eq.~\eqref{eq:adaptedguess}) and the deflated PCG with $2\times3$ vectors.}
    \label{fig:allsystems_simplified}
\end{figure*}

 \Cref{fig:allsystems_simplified} is the same as \Cref{fig:allsystems_compare}, but for the simplified setting. Here two triplets of the eigenvectors are constructed following the procedure described above.

\section{Ingredients of the proposed procedure}
\label{sec:ingredients}

We present in this section two ingredients of the proposed procedure in more detail. Namely, we discuss approaches for estimating the eigenpairs from the computed basis of the Krylov subspace and approaches for combining the deflation of the approximate eigenvectors with another preconditioner.  To facilitate presentation, we simplify the notation in this section.
 
\subsection{Approximating the eigenvalues using Krylov subspace methods}
\label{sec:eigvalueapprox}
We present first the Rayleigh--Ritz approximation, which is used in the Arnoldi and Lanczos algorithms to approximate the eigenvalues of a general nonsingular or, respectively, a hermitian matrix. Then, we recall the Arnoldi and Lanczos algorithms, and finally, we briefly comment on their restarted variants.

The methods discussed below do not represent an exhaustive overview of methods for approximating several eigenvalues and the associated eigenvectors. The omitted methods include the Jacobi--Davidson method \citep{SleVdV00}, for example, which proved to be particularly efficient for approximating the inner part of the spectrum. For a survey of the methods and a list of references, see \cite{Sor02}, for instance .

\subsubsection{Ritz value and harmonic Ritz value approximations}
\label{sec:Ritz}

For a subspace $\mathcal{S}\subset \mathbb{C}^n$, we call $y\in\mathcal{S}$ a \emph{\textup{Ritz vector}} of $A$ with \emph{\textup{Ritz value}} $\theta$ if
\begin{equation}
        Ay - \theta y \perp \mathcal{S}.
\end{equation}
When a (computed) basis $V_j$ of $\mathcal{S}$ is used and $y = V_j w$ is set, the above relation is equivalent to solving
\begin{equation}
\label{eq:RitzEP}
        V_j^\top A V_j w = \theta V_j^\top V_j w.
\end{equation}

Ritz values are known to approximate the extreme eigenvalues of $A$ well. When an approximation to the interior eigenvalues is required, it can be preferable to compute the harmonic Ritz values. The term \emph{\textup{harmonic Ritz values}} was introduced in~\cite{PaiPardVor95}, where references to previous works using this approximation can be found.
Following \cite{desturler06:recycleGMRES}, we define harmonic Ritz values as the Ritz values of $A^{-1}$ with respect to the space $A\mathcal{S}$,
\begin{equation}
        \widetilde{y} \in A\mathcal{S}, \quad A^{-1}\widetilde{y} - \widetilde{\mu} \widetilde{y} \perp A\mathcal{S}.
\end{equation}
We call $\widetilde{\theta} \equiv 1/\widetilde{\mu}$ a \emph{\textup{harmonic Ritz value}} and $\widetilde{y}$ a \emph{\textup{harmonic Ritz vector}}. When $V_j$ is a basis of $\mathcal{S}$ and $\widetilde{y} = V_j \widetilde{w}$, this relation can be represented as
\begin{equation}
\label{eq:harmRitzEP}
        V_j^\top A^\top V_j \widetilde{w} = \widetilde{\mu} V_j^\top A^\top A V_j \widetilde{w} \ \  \Longleftrightarrow \ \ V_j^\top A^\top A V_j \widetilde{w} = \widetilde{\theta} V_j^\top A^\top V_j \widetilde{w}.
\end{equation}
For the properties of the harmonic Ritz value approximations and the relationship with the iteration polynomial in MINRES method, see~\cite{PaiPardVor95}.

\begin{remark} There are various presentations and definitions in the literature of the harmonic (Rayleigh–)Ritz procedure; it is often introduced to approximate eigenvalues close to a target $\tau\in\mathbb{C}$. For example, \cite{Wu17} prescribes the procedure by
\begin{equation}
\label{eq:harmRitzEP-Wu17}
        \widetilde{y} \in \mathcal{S}, \quad A \widetilde{y} - \widetilde{\theta} \widetilde{y} \perp (A -\tau I)\mathcal{S},
\end{equation}
where $I$ is the identity matrix. With $\widetilde{y} = V_j \widetilde{w}$, this corresponds to the generalized eigenvalue problem 
\begin{equation}
        V_j^\top (A -\tau I)^\top (A -\tau I) V_j \widetilde{w} = \big(\widetilde{\theta} - \tau\big) \big(V_j^\top (A -\tau I)^\top V_j \big) \widetilde{w},
\end{equation}
which becomes for $\tau = 0$ exactly the right-hand side equality in Eq.~\eqref{eq:harmRitzEP}.
\end{remark}

We note that harmonic Ritz approximation is also often used in the Krylov subspace recycling methods to approximate the smallest (in magnitude) eigenvalues and the associated eigenvectors.

Finally, we comment on the (harmonic) Ritz approximation in the case when we wish to compute the eigenvalues of the matrix~$A$ preconditioned from the left by~$M$. 
In the general case, when only $A$ and $M$ are assumed to be nonsingular, the Ritz and harmonic Ritz approximation are applied as above by just replacing in the formulas~$A$ by $M^{-1}A$.
When the matrix~$A$ is hermitian and the preconditioner~$M$ is SPD, there is also another option. First, we note that the matrix $M^{-1}A$ is not hermitian, but is self-adjoint with respect to the inner product induced by~$M$, that is,
\begin{equation}
        (v, M^{-1}A w)_{M} = (M^{-1}A v, w)_{M}, \qquad \forall v, w,
\end{equation}
where $(v, w)_{M} \equiv v^\top M w$. This allows in the definition of Ritz and harmonic Ritz approximation replacing~$A$ by $M^{-1}A$ and the standard inner product by the inner product induced by the matrix~$M$, giving
\begin{equation}
        y\in\mathcal{S}, \quad M^{-1}Ay - \theta y \perp_M \mathcal{S}\end{equation}
or
\begin{equation}
    \widetilde{y} \in M^{-1}A\mathcal{S}, \quad (M^{-1}A)^{-1}\widetilde{y} - \widetilde{\mu} \widetilde{y} \perp_M M^{-1}A\mathcal{S}\end{equation}
respectively. The corresponding algebraic problems with $y = V_jw$ and  $\widetilde{y} = V_j\widetilde{w}$ are
\begin{equation}
        V_j^\top A V_j w = \theta V_j^\top M V_j w,
\end{equation}
and
\begin{equation}
    V_j^\top A^\top M^{-1} A V_j \widetilde{w} = (1/\widetilde{\mu}) V_j^\top A^\top V_j \widetilde{w},\end{equation}
respectively. The problems above involve hermitian matrices only.

\subsubsection{Arnoldi and Lanczos methods}
\label{sec:arnoldilanczos}

Arnoldi and Lanczos algorithms for approximating the eigenvalues of a general nonsingular or a hermitian matrix are based on a Ritz approximation with setting $\mathcal{S} = \mathcal{K}_j(A, v_1) = \mbox{span}(v_1, A v_1, \ldots, A^{j-1} v_1)$, the $j$th Krylov subspace. The methods compute an orthogonal basis~$V_j$ of $\mathcal{S}$ such that
\begin{equation}
        A V_j = V_j T_j + \beta v_{j+1} e_j^\top,
\end{equation}
where $e_j$ is the last column vector of the identity matrix (of size~$j$) and $V_j^\top V_j = I$, $V_j^\top v_{j+1} = 0$. Consequently, the eigenvalue problem in Eq.~\eqref{eq:RitzEP} corresponding to the Ritz approximation reads
\begin{equation}
        T_j w = \theta w.
\end{equation}
The matrix~$T_j$ is available during the iterations. The standard use of the Arnoldi and Lanczos method for eigenvalue approximation consists of solving the above problem and setting the pairs $(\theta, V_j w)$ as the computed approximations. 

The Ritz approximation can be replaced by the harmonic Ritz approximation. Then, the matrices in Eq.~\eqref{eq:harmRitzEP} become
\begin{equation}
        V_j^\top A^\top A V_j = T_j^\top T_j + \beta^2 e_j e_j^\top, \qquad V_j^\top A^\top V_j = T_j^\top.
\end{equation}

\begin{remark}
        The Lanczos algorithm is a variant of the Arnoldi algorithm for a hermitian $A$. The matrix~$T_j = V_j^\top A V_j$, which is in the Arnoldi method upper Hessenberg, is then also hermitian. Consequently, it is tridiagonal, which means that in each step of the Lanczos method, we orthogonalize the new vector only against the two previous vectors. This ensures that the computational cost of each iteration is fixed, and only when the eigenvalues are to be approximated, storing three vectors $v_{j-1}$, $v_{j}$ and $v_{j+1}$ is sufficient instead of handling the full matrix~$V_j$. The assumption on exact arithmetic is crucial here, however. In finite precision computations, the global orthogonality is typically quickly lost, which can cause several stability issues.
\end{remark}

As noted above, an orthonormal basis~$V_j$ of $\mathcal{S}$ is advantageous for the Ritz approximation. For the harmonic Ritz approximation applied to an SPD matrix~$A$, an \mbox{$A$-ortho}\-nor\-mal basis can instead be constructed, which ensures that the matrix $V_j^\top A^\top V_j = V_j^\top A V_j$ on the right-hand side of Eq.~\eqref{eq:harmRitzEP} is equal to the identity. An $A$-orthonormal basis of a Krylov subspace can be constructed within the iterations of conjugate gradient method using the search direction vectors.

\medskip
The Arnoldi method can also be applied to the preconditioned matrix $M^{-1}A$ to compute an orthonormal basis~$V_j$ of the associated Krylov subspace $\mathcal{K}_j(M^{-1}A, M^{-1}v_1)$, giving\begin{equation}
        M^{-1} A V_j = V_j T_j + \beta v_{j+1} e_j^\top, \quad V_j^\top V_j = I, \quad V_j^\top v_{j+1} = 0.
\end{equation}
For a hermitian~$A$ and an SPD preconditioner~$M$, we can apply the Lanczos method following the comment made above, using the matrix~$M^{-1}A$ and the inner product $(\cdot, \cdot)_M$ induced by~$M$ instead of the standard euclidean  $(\cdot, \cdot)$, giving
\begin{equation}
        M^{-1} A V_j = V_j T_j + \beta v_{j+1} e_j^\top, \quad V_j^\top M V_j = I, \quad V_j^\top M v_{j+1} = 0.
\end{equation}
The computed basis~$V_j$ is therefore $M$-orthonormal.

\subsubsection{Restarted variants}
\label{sec:restartedvariants}

The number of iterations necessary to converge is not a~priori known in Arnoldi and Lanczos algorithms, and it can in general be very high. High iteration counts require a large memory to store the basis vectors, and when a full reorthogonalization is used, a high computational effort because of the growing cost of the reorthogonalization in each step. The idea behind implicitly restarted variants is to limit the dimension of the search space~$\mathcal{S}$. This means that the iterations are stopped after a (fixed) number of steps, the dimension of the search space is reduced while maintaining its (Krylov) structure, and the Arnoldi/Lanczos iterations are resumed.

Several restarted variants are described in the literature (a detailed description is beyond the scope of this paper, however): the implicitly restarted Arnoldi (IRA, \cite{Sor92}), the implicitly restarted Lanczos (IRL, \cite{CalReiSor94}), or the Krylov--Schur method (\cite{Ste01,WuSim00}).

The estimation of the spectrum of~$A$ is possible within the GMRES, MINRES, and CG iterations (applied to solve a system with~$A$) because they are based on Arnoldi (GMRES and MINRES) or Lanczos (CG) algorithms. In contrast, a combination of restarted variants with solving a linear algebraic system is, to the best of our knowledge, not described in the literature.

\subsection{Deflation and two-level preconditioners}
\label{sec:deflationvariants}

In this section we first discuss a deflation preconditioner for Krylov subspace methods that can be regarded as eliminating the effect of several (given) vectors from the operator or, equivalently, augmenting by these vectors the space in which we search for an approximation. Then we describe a combination of the deflation preconditioner with another preconditioner that is commonly used in practice.

\smallskip
The Krylov subspace methods (in particular CG, \cite{CG}, and GMRES, \cite{GMRES}) are well-known for their minimization (optimal) properties over the consecutively built \emph{\textup{Krylov subspace,}}
\begin{equation}
        \mathcal{K}_j(A, v) = \mbox{span}\{v, Av, A^2 v, \ldots, A^{j-1} v\}.
\end{equation}
A question then arises: given some other subspace $\mathcal{U}$, can we modify the methods such that they have the same optimal properties over the union of $\mathcal{K}_j(A, v)$ and $\mathcal{U}$ (which is often called an \emph{\textup{augmented}} Krylov subspace)?
The answer is positive and the implementation differs according to the method: it is straightforward for GMRES and requires more attention for CG. Hereafter, we denote by $I$ the identity matrix and by $Z$ the basis of $\mathcal{U}$.

The deflation in GMRES method is often (see, e.g., GCROT by \cite{morgan95:GMRES-DR}) considered as a remedy to overcome the difficulties caused by restarts: for computational and memory restrictions, only a fixed number of GMRES iterations is typically performed, giving an approximation that is then used as the initial vector for a new GMRES run. In GCROT, several vectors are saved and used to augment the Krylov subspace built after the restart. The GMRES method with the deflation was used to solve a sequence of linear algebraic systems in~\cite{desturler06:recycleGMRES}, for example.

The augmentation of the Krylov subspace in CG is more delicate because the original CG method can only be applied to an SPD matrix. The first such algorithm was proposed in~\cite{nicolaides87:deflated-cg} and~\cite{Dos88}. We note that it includes the construction of the conjugate projector
\begin{equation}
        P_{c.proj.} = Z(Z^\top A Z)^{-1} Z^\top A
,\end{equation}
and in each iteration, the computation of the preconditioned search direction $q_i = (I-P_{c.proj.})\,p_i$ and of the vector $A q_i$. 
The latter can be avoided at the price of storing $Z$ \emph{and}~$A Z$ and performing additional multiplication with~$A Z$. In both variants, the cost of a single iteration is higher than the cost of one standard CG iteration.

The combination of a preconditioner with a deflation is widely studied in the literature and therefore we present this only briefly; more details and an extensive list of references can be found in the review paper by~\cite{vuik:deflation09}, for instance. The preconditioner stemming from the combination of a (typically relatively simple) traditional preconditioner with the deflation is called a \emph{\textup{two-level preconditioner}}. As shown in~\cite{vuik:deflation09}, this shows an analogy with multilevel (multigrid) and domain decomposition methods. While the traditional preconditioner aims at removing the effect of the largest (in magnitude) eigenvalues, the deflation (projection-type preconditioner) is intended to remove the effect of the smallest eigenvalues. Common choices for the traditional preconditioner are block Jacobi, (restricted) additive Schwarz method, and incomplete LU or Cholesky factorizations. In many applications, two-level preconditioners proved to be efficient in the CMB data analysis (see, e.g., \cite{GriStoSzy12,SzyGriSto14}).

We now present the combination of the traditional and projection-type (deflation) preconditioners following the discussion and notation of~\cite{vuik:deflation09}. Hereafter, we assume that the system matrix~$A$ and the traditional preconditioner~$M$ are SPD.
We note that some of the pre\-con\-di\-tio\-ners~$\mathcal{P}_{\Box}$ mentioned below are not symmetric. However, their properties allow us to use them (with possible modification of the initial vector) as left preconditioners in PCG; see~\cite{vuik:deflation09} for details.

Let the deflation space span the columns of the matrix~$Z$. We denote
\begin{equation}
\label{eq:defprec}
        P \equiv I - A Q, \qquad Q \equiv Z\big(Z^\top A Z\big)^{-1}Z^\top.
\end{equation}
Two-level preconditioners based on the deflation are given as
\begin{equation}
        \mathcal{P}_{\text{DEF1}} \equiv M^{-1}P, \qquad \mathcal{P}_{\text{DEF2}} \equiv P^\top M^{-1}.
\end{equation}
Other preconditioners can be determined using the \emph{\textup{additive}} combination of two (SPD) preconditioners $C_1$, $C_2$ as
\begin{equation}
        \mathcal{P}_{\text{add}} \equiv C_1 + C_2,
\end{equation}
or, using the \emph{\textup{multiplicative}} combination of the preconditioners, as
\begin{equation}
        \mathcal{P}_{\text{mult}} \equiv C_1 + C_2 - C_2 A C_1.
\end{equation}
Three variants of two-level preconditioners are derived by choosing aan dditive or multiplicative combination and setting $C_1 = M^{-1}$, $C_2 = Q$, or $C_1 = Q$, $C_2 = M^{-1}$. Other preconditioners can be derived using the multiplicative combination of three SPD matrices (see \cite[~]{vuik:deflation09}.

The variants of the two-level preconditioner mentioned above differ in the implementation cost and also in the numerical stability; see~\cite{vuik:deflation09}.
The variant~$\mathcal{P}_{\text{DEF1}}$, which is often used in the procedures for solving the sequences of linear systems (see, e.g., \cite{saad00:deflated-cg}), was found to be cheap but less robust, especially with respect to the accuracy of solving the coarse problem with the matrix $Z^\top A Z$ and with respect to the demanded accuracy. The conclusion drawn in~\cite{vuik:deflation09} is that "\mbox{A-DEF2} seems to be the best and most robust method, considering the theory, numerical experiments, and the computational cost". Therefore the preconditioner $\mathcal{P}_{\text{A-DEF2}}$,
\begin{equation}
\label{eq:adef2}
        \mathcal{P}_{\text{A-DEF2}} \equiv M^{-1} + Q - Q A M^{-1} = P^\top M^{-1} + Q,
\end{equation}
is of interest, in particular in the cases where the dimension of the deflation space (equal to the number of columns in~$Z$) is high and/or the matrix~$M^{-1}A$ is ill-conditioned.

As noted in~\cite{saad00:deflated-cg}, the gradual loss of orthogonality of the computed residuals with respect to the columns of~$Z$ can cause stagnation, divergence, or erratic behaviour of errors within the iterations (see also the comment in \cite[~]{vuik:deflation09}. The suggested remedy in this case consists of reorthogonalizing the computed residuals as
\begin{equation}
\label{eq:reortores}
        r_j := W r_j, \qquad W \equiv I - Z(Z^\top Z)^{-1} Z^\top.
\end{equation}
However, no such instabilities were observed in our experiments, and the results depicted throughout the paper are for the standard (non-reorthogonalized) variant.

\section{Full algorithm}
\label{sec:fullalg}

In this section we provide the pseudo-codes for the deflated PCG (\Cref{alg:deflPCG}), the subspace recycling (\Cref{alg:subrec}), and for the full procedure (\Cref{alg:thealg}) proposed in this paper and used in the numerical experiments in \Cref{sec:numexp}.

\begin{algorithm}[H]
\caption{deflated PCG (variant "def1")\label{alg:deflPCG}}
\begin{algorithmic}
\setstretch{1.2}
\Function{deflPCG}{$\A$, $\M$, $\sigg{b}$, $\sigg{s}^{(0)}$, $Z$, $\dim_p$, $j_{\text{max}}$}
\State $Q = Z(Z^\top \A Z)^{-1}Z^\top$ \Comment{in practice, we save $\A Z$ to use later}
\State $\sigg{r}^{(0)} = (I-\A Q)(\sigg{b} - \A \sigg{s}^{(0)})$  \Comment{with saved $\A Z$, $Q\A$ and $\A Q$ can be computed without applying $\A$}
\State $\sigg{p}^{(0)} = \sigg{r}^{(0)}$
\State $\widetilde{\sigg{r}}^{(0)} = \M^{-1} \sigg{r}^{(0)}$
\For{$j = 0, \dots, j_{\text{max}}$}
  \State $\sigg{w}^{(j)} = (I-\A Q)(\A \sigg{p}^{(j)})$
  \If{$j \leq \dim_p$} 
    \State \textbf{save} $(I-Q\A) \sigg{p}^{(j)}$ into $\widetilde{Z}$
    \Comment{in practice, we also save  $\sigg{w}^{(j)}$ to avoid computing $\A\widetilde{Z}$ later}
  \EndIf
  \State $\gamma^{(j)} = {(\widetilde{\sigg{r}}^{(j)},\sigg{r}^{(j)})}/{(\sigg{p}^{(j)}, \sigg{w}^{(j)})}$
  \State $\sigg{s}^{(j+1)} = \sigg{s}^{(j)} + \gamma^{(j)} \sigg{p}^{(j)}$
  \State $\sigg{r}^{(j+1)} = \sigg{r}^{(j)} - \gamma^{(j)} \sigg{w}^{(j)}$
  \State \textbf{check} the stopping criterion
  \State $\widetilde{\sigg{r}}^{(j+1)} = \M^{-1} \sigg{r}^{(j+1)}$
  \State $\delta^{(j)} = {(\widetilde{\sigg{r}}^{(j+1)},\sigg{r}^{(j+1)})}/{(\widetilde{\sigg{r}}^{(j)},\sigg{r}^{(j)})}$
  \State $\sigg{p}^{(j+1)} = \widetilde{\sigg{r}}^{(j+1)} + \delta^{(j)} \sigg{p}^{(j)}$
\EndFor
\State $\sigg{s}^{(final)} := Q\sigg{b} + (I-Q\A) \sigg{s}^{(j)}$
\State 
\Return $\sigg{s}^{(final)}$; $\widetilde{Z}$ 
%\State \Comment{we return the computed approximation and the saved search vectors}
\State \Comment{to be efficient, we also return $\A Z$ and the vectors $\{\sigg{w}^{(j)}\}$}
\EndFunction
\end{algorithmic}
\end{algorithm}

\begin{algorithm}[H]
\caption{subspace recycling (variant "ritz")\label{alg:subrec}}
\begin{algorithmic}
\setstretch{1.2}
\Function{SubspRec}{$\A$, $\M$, $U$, $k$}
\State $F = U^\top \M U$
\State $G = U^\top \A U$
\Comment{in practice, we reuse $\A U$ saved during the deflated PCG}
\State \textbf{solve} the generalized eigenvalue problem $GY = \mbox{diag}(\lambda_i)FY$
\State \textbf{take} $k$ smallest $\lambda_i$ and the respective columns of $Y$, $Y_k$
\State
\Return $Z = UY_k$
\EndFunction
\end{algorithmic}
\end{algorithm}

\begin{algorithm}[H]
\caption{full algorithm of the procedure\label{alg:thealg}}
\begin{algorithmic}
\setstretch{1.4}
\Require $\beta_0$, $\sigg{s}^{(0)}_{\beta_0}$
\Require $k$, $\dim_p$
\State set $Z := []$
\For{$i = 0, \dots, i_{\text{max}}$}
    \State \textbf{assembly} the system matrix $\A$, right-hand side $\sigg{b}$, and the preconditioner $\M$ corresponding to $\beta_i$
    %\State \textbf{if} $i>0$ \textbf{then} $s^{(0)}_{\beta_i} = (\widetilde{M}_{\beta_i})^\dagger \widetilde{M}_{\beta_{i-1}} s^{(final)}_{\beta_{i-1}}$
    \State \textbf{if} $i>0$ \textbf{then} $\sigg{s}^{(0)}_{\beta_i} = (\widetilde{M}_{\beta_i})^\dagger \sigg{s}^{(0)}_{\beta_i}$
    %\State run deflated PCG 
    \State \Call{deflPCG}{$\A$, $\M$, $\sigg{b}$, $\sigg{s}^{(0)}_{\beta_i}$, $Z$, $\dim_p$, $j_{\text{max}}$} $\longrightarrow$ ($\sigg{s}^{(final)}_{\beta_{i}}$, $\widetilde{Z}$) 
    \State \textbf{check} the stopping criterion for $\beta_i$, \textbf{exit} if $i=i_{\text{max}}$
    \State $\sigg{s}^{(0)}_{\beta_{i+1}} = \widetilde{M}_{\beta_{i}} \sigg{s}^{(final)}_{\beta_{i}}$
    \State ({determine} $\beta_{i+1}$) \Comment{ considered here as a black box}
    %\State run subspace recycling to determine new $Z$
    \State \Call{SubspRec}{$\A$, $\M$, $U = [Z, \widetilde{Z}]$, $k$} $\longrightarrow$ $Z$
\EndFor
\end{algorithmic}
\end{algorithm}

\section{Results for an alternative sequence of mixing parameters from a Monte Carlo sampling}
\label{sec:DominicsBetas}

\begin{figure}[ht]
    \centering
    \includegraphics[width=0.7\linewidth]{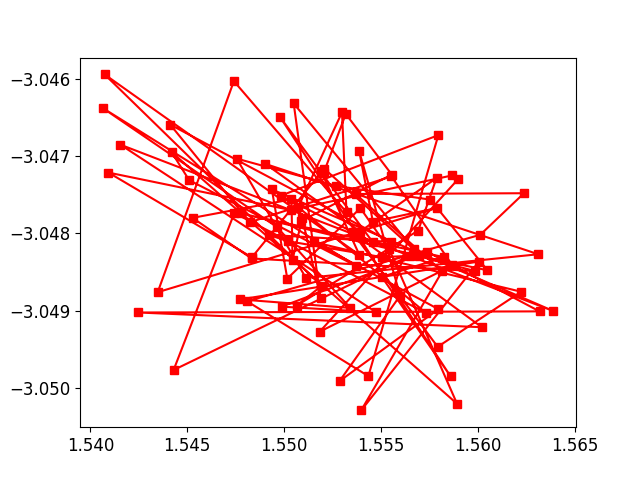}
    \caption{Plot of a sequence of the spectral parameters~$\beta_i = [\beta_{i,s}, \beta_{i,d}]$ drawn through a Monte Carlo sampling technique and used as an alternative test case in the numerical experiments described in~\Cref{sec:DominicsBetas}.}
    \label{fig:betasDominic}
\end{figure}

In this section we provide results for a sequence of spectral parameters generated by a Monte Carlo sampling algorithm. In realistic circumstances, these sequences may contain up to many thousand samples, but for computational efficiency, we here restrict ourselves to a subsequence made of merely $30$ elements. We use them in order to demonstrate the performance of the proposed approach on a sequence with realistic properties but sufficiently different than those of the sequences encountered in the likelihood maximization process. We emphasize that it is not yet our purpose here to validate the method on the actual application, and more work is needed to achieve this goal, see Sect.~\ref{sec:conclusion}. 
The sequence is depicted in \Cref{fig:betasDominic}. It was produced using the publicly available software {\sc fgbuster}\footnote{{\sc fgbuster}: {\tt https://github.com/fgbuster}} and indeed shows qualitatively a very different behavior than that of our standard case displayed in Fig.~\ref{fig:betas}.

\begin{figure*}[ht]
    \centering
    \includegraphics[width=0.9\textwidth]{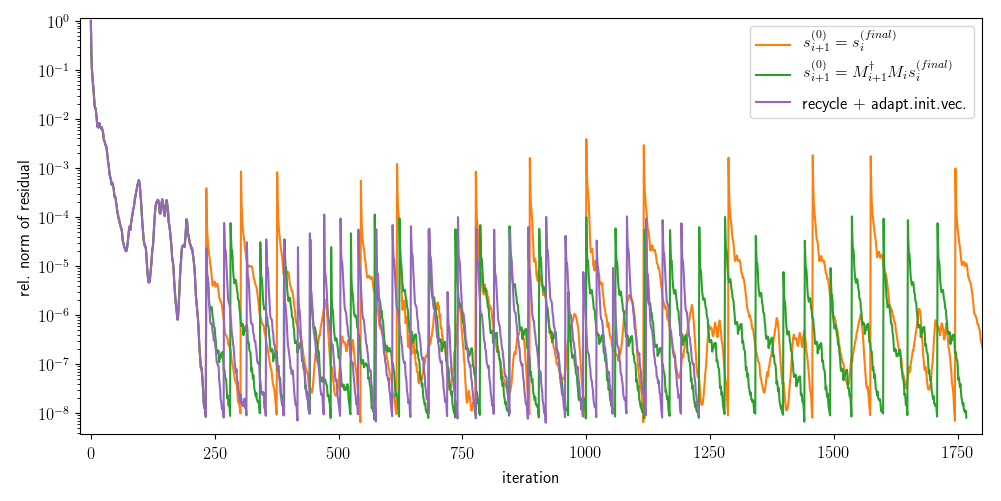}
    \caption{Comparison of the PCG with different choices of initial guess (as in Eqs. \eqref{eq:continuation} and~\eqref{eq:adaptedguess}) and the PCG with the subspace recycling (together with the choice of the initial guess as in Eq.~\eqref{eq:adaptedguess}). For the recycling, we consider $k=10$ eigenvectors approximated using $\dim_p = 100$. The convergence for the whole sequence when the initial guess is as in Eq.~\eqref{eq:continuation} (the yellow line) requires 4010 iterations.}
    \label{fig:allsystems_compareDominic}
\end{figure*}

In \Cref{fig:allsystems_compareDominic} and in \Cref{tab:allsystems_compareDominic}, we compare the results obtained in this case by applying the various techniques discussed and proposed in this work.

\begin{table}[htp]
\centering
\begin{tabular}{l|ccc}
 & \multicolumn{3}{c}{\#MatVecs} \\ 
 & iteration & deflation & total\\ \hline
$\sigg{s}^{(0)}_{\beta_{i+1}}$ as in \eqref{eq:continuation} & 4010 & 0 & 4010\\
$\sigg{s}^{(0)}_{\beta_{i+1}}$ as in \eqref{eq:adaptedguess} & 1768 & 0 & 1768\\
recycle + $\sigg{s}^{(0)}_{\beta_{i+1}}$ as in \eqref{eq:adaptedguess} & 1228 & 290 & 1518
\end{tabular}
\caption{Number of matrix-vector products (MatVecs) for different techniques as in \Cref{fig:allsystems_compareDominic}.}
    \label{tab:allsystems_compareDominic}
\end{table}

\bibliographystyle{aa}
\bibliography{./37687}

\end{document}